\documentclass[12pt]{iopart}
\usepackage{amssymb}
\usepackage{amsfonts}
\usepackage{graphics}
\usepackage{graphicx}
\usepackage{subfigure}
\usepackage{multirow}
\usepackage{color}
\usepackage{textcomp}
\usepackage{float}
\newcommand{\lt}{\left}
\newcommand{\rt}{\right}

\newcommand{\cL}{\mathcal{L}}
\newcommand{\cA}{\mathcal{A}}

\newcommand{\cT}{\mathcal{T}}
\newcommand{\R}{\mathbb{R}}

\newcommand{\vecp}{{\mathrm{vec}}}
\usepackage{algorithm}
\usepackage{algorithmicx}
\usepackage[noend]{algpseudocode}
\newcommand*\Let[2]{\State #1 $\gets$ #2}
\algrenewcommand\algorithmicrequire{\textbf{Input:}}
\algrenewcommand\algorithmicensure{\textbf{Output:}}

\usepackage{iopams}  
  \expandafter\let\csname equation*\endcsname\relax
  \expandafter\let\csname endequation*\endcsname\relax
\usepackage{amsmath}
\begin{document}

\title[Efficient travel-time tomography]{Relaxation algorithms for matrix completion, 
with applications to seismic travel-time data interpolation.}

\author{Robert Baraldi$^{1}$, Carl Ulberg$^{2}$, Rajiv Kumar$^{3}$, Kenneth Creager$^{2}$ and Aleksandr Aravkin$^{1}$}

\address{$^{1}$Department of Applied Mathematics, University of Washington\\
$^{2}$Department of Earth and Space Sciences, University of Washington\\
$^{3}$School of Earth and Atmospheric Sciences, Georgia Institute of Technology\\
}
%\ead{submissions@iop.org}
%\vspace{10pt}
%\begin{indented}
%\item[]August 2017
%\end{indented}

\begin{abstract}
Travel time tomography is used to infer the underlying three-dimensional wavespeed structure of the Earth 
by fitting seismic travel time data collected at surface stations. 
Data interpolation and denoising techniques are important pre-processing steps that 
use prior knowledge about the data, including parsimony in the frequency and wavelet domains, 
low-rank  structure of matricizations, and local smoothness. 
  
We show how local smoothness structure can be combined with low rank  
constraints using level-set optimization formulations, and develop a new relaxation algorithm 
that can efficiently solve these joint problems. 
In the seismology setting, we use the approach 
to interpolate missing stations and de-noise observed stations. 
The new approach is competitive with alternative algorithms,
and offers new functionality to interpolate observed data using both smoothness and low rank structure in the presence of data fitting constraints.  

%We compare this to approach to available formulations and algorithms, and show that the new approach is competitive and has expanded functionality, in particular for non-smooth and non-convex constraints. 
\end{abstract}
 
\submitto{\IP}
% \maketitle
% \ioptwocol
%

\section{Introduction}
% \textbf{Travel-time tomography in seismic inverse problems}. \\
Travel-time tomography is used to determine the underlying structure of the earth and how seismic waves propagate through that structure. This weakly nonlinear problem is formulated as a data-fitting inverse problem 
and solved using iterative optimization techniques. Data quality and availability are key constraints and can drastically influence the merit of the results~\cite{phillips1991comparison,landro2002uncertainties,landro2008effect}. 
Hence, researchers often use prior information to denoise and interpolate the data prior to inversion. Parsimonious representations~\cite{candes2006near,Recht2010guaranteed} 
of the data in transform domains such as Fourier~\cite{sacchi} and Curvelet~\cite{herrmann2008non} have been used in exploration seismology~\cite{hennenfent2008simply,kumar2015efficient}, along with low-rank representations~\cite{aravkin2014fast,splittingschemes}. 

We focus on regional seismology --- in particular, we want to analyze data obtained by the Imaging Magma Under St. Helens (iMUSH) project,  a multi-year effort to image and infer the architecture of the greater Mount St. Helens, WA, magmatic system \cite{2017AGUFMcarl,2017AGUkiser}. 
The data collected from this experiment comprises {active experiments}~\cite{kiser2016vpvs} that initiate seismic tremors and record their waves/reverberations at different stations, 
and {passive experiments} that record seismic waves from natural earthquakes. Collected data has a range of fidelities, and different subsets inform different parameters for 3D $S$- and $P$-wave velocities and hence geometries of the structure underlying Mount St. Helens. %inversion methods 
%use \texttt{struct3DP} (provide some citation for solver) to 
Travel times are inverted to obtain 3D seismic velocity models %using a 3D eikonal equation solver
and image complex subsurface structures, including low velocity zones.
The efficacy of the approach depends on the quality of the data, which is affected by noise from roads, streams, ocean waves, and wind. % and more that may affect the seismograms. 
%The system is very large, with millions of parameters, and extremely sparse across the domain. However, it is well known that
The signal to noise ratio (SNR) of the data depends on the earthquake size, distance from the seismometer, and attenuation structure of the earth.  
Data is scarce because many landscape features are impassable; typical station spacing is 10km. 
Arrival times are chosen by operators, who assign uncertainties for particular earthquakes based on confidence in seismic readings.
Low-magnitude events in particular are often difficult to distinguish from noise. 
        % \subsection*{Discuss ways of overcoming these limitations and why they might not work} 
%Installing more sensors may overcome data sparsity problems as well as some noise issues; however, the issue of impassable terrain remains, as well as sensor placement. In addition, the iMUSH project has been completed, so no more sensors can be installed or used (outside of the permanent ones). Running the inverse problem at a higher resolution may resolve some resolution issues, but that does not necessarily mean results are correct in some locations. Usually a checkerboard test is done to see where exactly the model has good resolution versus bad resolution by examining how sensitive it is to perturbations, but this does not help in resolving the areas of low-resolution. Running the model may also result in needless computational expense; \texttt{struct3DP} takes several days (?) to run already. 

Our goal is to denoise available data, and spatially interpolate to additional locations 
that might increase resolution and stability of the overall inverse problem. 
%The goal of this paper is to show that a new low-rank matrix completion formulation can accurately reproduce seismic data through a synthesis of local and global information. 
    % \subsection{Where we are going} 
%Low-rank matrix completion and denoising is attractive because it allows the enunciation of complex data features with a few components. For example, \cite{mallat2008wavelet,hennenfent2010GEOPnct} show that wavelet and curvelet transforms can be expressed with few principal components, thereby lending missing/noisy data to be re-expressed with linear combinations of known/trusted data. In addition to more famous examples such as the Netflix problem, 
%The matricization of data naturally leads to low-rank interpolation schemes \cite{aravkin2014fast,aravkin2014fast,splittingschemes}, as storing data in matrix form phrases it in the context of a linear basis. 
Since we know the data is corrupted by uncertainty and noise,
% of the noise-floor (through uncertainties provided by operators), 
the problem is a good match for {level-set} optimization formulations~\cite{aravkin2016level} 
that minimize a regularizer subject to a prescribed level of data fit. %measured in a particular norm. 
%We also allow both local smoothness regularization and low-rank regularization. 
We propose a relaxation formulation that allows misfit constraints while
(1) penalizing rank of a tessellation that captures redundancy of features across sources, 
and (2) enforcing smooth features consistent with the underlying physical model. 
We develop an efficient block-coordinate algorithm for this formulation, and compare the approach  
with a variety of competing formulations and algorithms. 

 The paper proceeds as follows. 
  In Section~\ref{sec:Review} we review relevant formulations and algorithms for interpolation and denoising. In Section~\ref{sec:Relax}, 
  we develop the extended model formulation and the relaxation method  to solve it. Section~\ref{sec:data} describes the data used and how it fits into a low-rank interpolation scheme. In Section~\ref{sec:results} we 
  evaluate the approach and compare it against alternatives for denoising and interpolating Mount St. Helens data.

\section{Preliminaries}
\label{sec:Review}

In this section, we set up the notation, review prior information used in interpolation and denoising, 
and discuss different types of standard optimization formulations needed to implement such approaches.  

\subsection{Notation}
We use the following notation conventions in the paper. 
Lowercase variables ($x$) denote vectors, while uppercase variables ($X$) represent matrices. 
Calligraphic uppercase letters ($\cA$) are used for operators or functionals. 
The terms $(R_x, R_y)$ represent a receiver coordinate grid while $(S_x, S_y)$ represent the source coordinate grid. 
Single source/receiver coordinates are specified with an index $i$, as in $(S_{x_i}S_{y_i})$.
 % However, note that when referring to a specific source, it is written as $S_{x_i}S_{y_i}$. 
The variable $\Omega$ is used to represent the entire 4-dimensional source/receiver space and generically the interpolation space.

\subsection{Formulations for Low-rank Interpolation}
The goal of low-rank matrix completion is to accurately estimate the missing data entries of a matrix $X \in \R^{m\times n}$ from observed entries. Low-rank structure is often inferred if completed $X$ has few non-zero singular values and experiences entry repetition. Borrowing some formulation from \cite{aravkin2014fast}, let $\cT \subset \{1,\dots,n\} \times \{1, \dots, m\}$ be the set of observed entries, which is created by the sampling operator $\cA:\R^{n\times m} \rightarrow \R^{n\times m}$ via element-wise selection
\[
    \cA(X) = \begin{cases} X_{ij}, \quad \mbox{if} \ (i,j)\in \cT, \\ 0, \quad (i,j)\not\in\cT \end{cases}
\]
for which the observed data $b$ can be written as $b = \cA(X) + \epsilon$ for $\epsilon = \cA(\epsilon)\in \R^{n\times m}$. Here, $\epsilon$ represents the data corruption of observed entries via some noise distribution. 
In formulating the problem, the nuclear norm 
\[
\|X\|_* = \sum_{j =1}^{\min(n,m)} \hat{\sigma}_j (X),
\]
with $\hat\sigma_j$ the singular values of $X$, is used as a proxy for rank. The classic formulations that balance 
data fit with regularization are 
\begin{eqnarray}
\label{eq:tikhonov}
     &\min_{X\in\R^{n\times m}} \|X\|_* + \frac{1}{\sigma}\|\cA(X) - b\|_2\\
     \label{eq:ivanov}
         &\min_{X\in\R^{n\times m}}\quad \|\cA(X) - b \|_2 \quad \mbox{s.t.} \quad  \|X\|_* \leq \tau \\     
     \label{eq:morozov}
         &\boxed{\min_{X\in\R^{n\times m}} \|X\|_* \quad \mbox{s.t.}\quad \|\cA(X) - b \|_2 \leq \sigma}. 
\end{eqnarray}
These formulations are known as Tichonoff, Ivanov, and Morozov regularization~\cite{oneto2016tikhonov}, respectively. 
The misfit-constraint Morozov variant~\eqref{eq:morozov} is best suited for situations where a good estimate 
of the uncertainty, $\sigma$, is available. To simplify exposition, we will focus on Morozov-type formulations.

\subsection{Smoothness constraints}

In travel time tomography, smoothness and continuity between gridpoints is a reasonable prior, 
since the geological structure of the crust exhibits the traits of approximately homogeneous media. 
Smoothness is enforced by introducing a penalty term 
$\|\cL(X)\|_2^2$, with $\cL$ the discretization of the Laplacian operator. In 1D, it is a tridiagonal difference matrix operating on the vectorized $X$; 
in 2D, it has four off-diagonal elements.

Enforcing smoothness and low-rank structure combines local and global information. 
A  Morozov formulation, with $\gamma$ balancing the regularizers, is given by 
\begin{equation}
\label{eq:MorozovS}
    \min_X \|X\|_* + \frac{1}{2\gamma}\|\cL(X)\|_2^2 \quad \mbox{s.t.}\quad \ \|\cA(X) - b\|_2\leq \sigma.
\end{equation}

%Rephrasing this in $LR^T$ form yields
%\begin{eqnarray*}
%    &\min_{L, R} \frac12 \lt\|\begin{array}{c} L\\  R\end{array}\rt\|_F^2 + \frac{1}{2\gamma}\|\cL(LR^T)\|_2^2\\
%&\mbox{s.t.}\quad\ \|\cA(LR^T) - b\|_2\leq \sigma
%\end{eqnarray*}
%and employing the variable relaxation technique yields
%\begin{eqnarray}
%    &\min_{L, R,W} \frac12 \lt\|\begin{array}{c} L\\  R\end{array}\rt\|_F^2 + \frac{1}{2\gamma}\|\cL(W)\|_2^2 + \frac{1}{2\eta}\|W - LR^T\|_2^2\nonumber\\
%&\mbox{s.t.}\quad\ \|\cA(W) - b\|_2\leq \sigma.
%\label{eq:smoothconstrained}
%\end{eqnarray}
%

\subsection{Factorized Formulations}
\label{sec:factform}
Theoretical properties of matrix completion via nuclear-norm minimization have been extensively studied \cite{candes2010Matcompnoise,Candes2009matcompconvex}. 
The theoretical appeal of convex formulations is tempered by computational considerations --- algorithms that optimize $\|X\|_*$ require 
a full matrix decision variable, and full or partial singular value decompositions (SVDs)  
at each iteration. An efficient alternative is to use matrix factorization formulations~\cite{aravkin2014fast}, 
writing $X = LR^T$, with $L\in\R^{n\times k}, R\in \R^{m\times k}$. From \cite{Recht2010guaranteed}, 
we have the characterization 
\[
    \|X\|_* = \inf_{L, R: X = LR^T} \frac12 (\|L \|_F^2 + \|R\|_F^2), %= \inf_{L, R: X = LR^T} \frac12 \lt\| \begin{array}{c} L\\  R\end{array} \rt\|_F^2
\]
which allows us to replace $\|X\|_*$ in any formulation by $\frac12 (\|L \|_F^2 + \|R\|_F^2)$ for any factorization $X=LR^T$. 
For example, the three formulations~\eqref{eq:tikhonov}-\eqref{eq:morozov} become
\begin{eqnarray}
\label{eq:tikhonovF}
     &\min_{L\in\R^{n\times k},\ R\in\R^{m\times k}} \frac12  \|L\|_F^2 + \frac12 \|R\|_F^2 + \frac{1}{\sigma}\|\cA(LR^T) - b\|_2\\
     \label{eq:ivanovF}
         &\min_{L\in\R^{n\times k},\ R\in\R^{m\times k}}\quad \|\cA(LR^T) - b \|^2 \quad \mbox{s.t.} \quad  \|L\|_F^2 +  \|R\|_F^2 \leq 2\tau \\     
     \label{eq:morozovF}
         &\boxed{\min_{L\in\R^{n\times k},\ R\in\R^{m\times k}} \frac12  \|L\|_F^2 + \frac12 \|R\|_F^2  \quad \mbox{s.t.}\quad \|\cA(LR^T) - b \|_2 \leq \sigma},
\end{eqnarray}
where $k \ll \min(n, m)$, and the memory requirements are reduced from $mn$ to $k(n+m)$. 
No SVDs are required; formulation~\eqref{eq:tikhonovF} is smooth, formulation~\eqref{eq:ivanovF} requires 
simple projections onto the Frobenius-norm ball, and formulation~\eqref{eq:morozovF}
can be solved using~\eqref{eq:ivanovF} via root-finding as described by~\cite{aravkin2014fast}.

%This formulation possessess several advantages over the previous, unrelaxed formulation. The first is the easy separability between $L$, $R$, and $W$ updates. Notice that in \ref{eq:smoothconstrained}, the $L$ and $R$ updates are as they were before, as opposed to the unrelaxed equation which represents a difficult to differentiate and computationally inefficient method of minimization via L-BFGS. When solving the unconstrained version with this technique, the algorithm has to form and store a Hessian for $L$ and $R$ each iterate, thereby using an excess of computational storage. There is also a easily found solutions for $W$; as it can no longer be split into observed and non-observed points due to the regularization term, we can solve for $W$ directly if $\sigma = 0$ and use a root-finding algorithm for $\sigma>0$. Updating $W$ requires solving a trust-region subproblem with the regularized Laplacian as the linear system.

Our goal here is to solve the factorized Morozov formulation corresponding to~\eqref{eq:MorozovS}:
\begin{equation}
\label{eq:MorozovSF}
    \min_{L, R}\frac12  \|L\|_F^2 + \frac12 \|R\|_F^2 + \frac{1}{2\gamma}\|\cL(LR^T)\|_2^2\quad  \mbox{s.t.}\quad\  \| \cA(LR^T) - b\|_2\leq \sigma.
\end{equation}

In particular, this formulation incorporates both local and global structure, and gives a misfit target $\sigma$. It requires a new algorithm, since~\eqref{eq:MorozovSF} cannot be solved by the level-set approach of~\cite{aravkin2014fast}.  The only available alternative is to use Tichonoff-type formulations (see Section~\ref{sec:alt}). Our main technical contribution is a nonconvex splitting algorithm for problem~\eqref{eq:MorozovSF}, developed in the next section.

\section{Relaxed Joint Inversion}
\label{sec:Relax}

The main challenge of the factorized Morozov formulation~\eqref{eq:MorozovSF} is the data-misfit constraint. To solve
the problem, we propose a {\it relaxation} following the ideas of~\cite{zheng2018fast}. In particular, 
we introduce an auxiliary variable $W \approx LR^T$:
\begin{equation}
\label{eq:MorozovSFR}
    \min_{L, R, W}\frac12  \|L\|_F^2 + \frac12 \|R\|_F^2+ \frac{1}{2\gamma}\|\cL(W)\|_2^2 + \frac{1}{2\eta} \|W - LR^T\|_F^2 \quad  \mbox{s.t.}\quad\ \|\cA(W) - b\|_2\leq \sigma.
\end{equation}
Problem~\eqref{eq:MorozovSFR} is a relaxation for problem~\eqref{eq:MorozovSF}, since $W$ approximates $X = LR^T$; in particular $\|W-LR^T\| = \mathcal{O}(\eta)$. The salient modeling features  of~\eqref{eq:MorozovSF} are still preserved. We can now design a simple block-coordinate descent algorithm by iteratively optimizing in each of $(L,R,W)$, 
detailed in Algorithm~\ref{alg:block-descent-smooth}.

\begin{algorithm}[H]
\caption{Block-Coordinate Descent for~\eqref{eq:MorozovSFR}.}
\label{alg:block-descent-smooth}
\begin{algorithmic}[1]
\State{\bfseries Input:} $W_0, L_0, R_0$
\State{Initialize: $k=0$.}
\While{not converged}
\Let{$L_{k+1}$}{$\left(I +\eta R_k^TR_k\right)^{-1}(\eta W_kR_k)$} {\Comment{ {\scriptsize Solves $\frac12  \|L\|_F^2 + \frac{1}{2\eta} \|W_k - LR_k^T\|_F^2 $}}}
\Let{$R_{k+1}$}{$(\eta W_k^T L_{k+1})\left(I + \eta L_{k+1}^TL_{k+1}\right)^{-1}$} \Comment{{\scriptsize Solves $ \frac12 \|R\|_F^2 +  \frac{1}{2\eta} \|W_k - L_{k+1}R^T\|_F^2$}}\vspace{.1in}
\Let{$W_{k+1}$}{$\arg\min_{W}\frac{1}{2\gamma}\|\cL(W)\|_2^2 +  \frac{1}{2\eta} \|W - L_{k+1}R_{k+1}^T\|_F^2 \quad \mbox{s.t.}\quad \|\cA(W) - b\|_2\leq \sigma  $}
\Let{$k$}{$k+1$}
\EndWhile
\State{\bfseries Output:} $W_k, L_k, R_k$
\end{algorithmic}
\end{algorithm}

Steps $4$ and $5$ of Algorithm~\ref{alg:block-descent-smooth} are simple least squares updates; each minimizes~\eqref{eq:MorozovSFR} in $L$ and $R$ respectively, with the remaining variables held fixed. 
Algorithm~\ref{alg:block-descent-smooth} converges to a stationary point of~\eqref{eq:MorozovSFR} by~\cite[Theorem 4.1]{tseng2001convergence}. In particular, $f(L,R,W)$ has a unique minimum in each coordinate block 
with the remaining blocks held fixed, which satisfies condition (c) of the theorem. The uniqueness of the minima are clear from the closed form solutions in steps 4 and 5, and from the strong convexity 
of the $W$ subproblem in step 6.  Step $6$ is solved using an efficient root-finding method described below. The equivalent penalized problem  with penalty parameter 
$\lambda$ is given by 
\begin{equation}
\label{eq:wlam}
w(\lambda) := \arg\min_{w} \frac{1}{2\gamma}\|Lw\|_2^2 +  \frac{1}{2\eta} \|w - d_{k+1}\|_F^2 + \frac{\lambda}{2} \|Aw - b\|^2,
\end{equation}
where $w  = \mbox{vec}(W)$, $L$ is a sparse matrix that encodes the action of the Laplacian on $w$, $A$ 
is a sparse mask that pulls out the entries of $w$ to compare with $b$, and $d_{k+1} = \mbox{vec}(L_{k+1}R_{k+1}^T)$. 
The root finding method obtains the smallest value of $\lambda$ satisfying $\|\cA w(\lambda) - b\|_2\leq \sigma$. 

Taking the gradient of the objective defining $w(\lambda)$ in~\eqref{eq:wlam} and setting it equal to $0$, we find an explicit formula
\begin{equation}
\label{eq:tr}
w(\lambda) = \left(\frac{1}{\gamma} (L^TL) + \frac{1}{\eta} I + \lambda A^TA\right)^{-1}\left( \frac{1}{\eta}d_{k+1} + \lambda A^Tb\right).
\end{equation}
We need only find the smallest $\lambda \geq 0$ so that $\|Aw(\lambda)-b\|_2 = \sigma$. 
The special case $\lambda = 0$ occurs when the constraint is satisfied at the least squares solution $w(0)$ in~\eqref{eq:tr}:
\[
\left\|A  \left(\frac{1}{\gamma} (L^TL) + \frac{1}{\eta} I \right)^{-1}\left( \frac{1}{\eta}d_{k+1} \right)-b\right\| \leq \sigma,
\]
In all other cases, we have  
\[
f(\lambda) :=  \sigma - \|Aw(\lambda)-b\|_2, \quad f'(\lambda) = -\frac{\langle A^TAw(\lambda)-b,  \nabla_\lambda w(\lambda) \rangle.}{\|Aw(\lambda)-b\|_2}
\]
We compute the quantity $\nabla_\lambda w(\lambda)$ required to evaluate $f'(\lambda)$ using the complex step method~\cite{martins2001connection}, instead of differentiating~\eqref{eq:tr} directly. The root-finding 
update is given by 
\[
\lambda^+ := \lambda - \frac{f(\lambda)}{f'(\lambda).}
\]
%From the formula~\eqref{eq:tr}, Step $6$ of Algorithm~\ref{alg:block-descent-smooth} is equivalent to the quadratic trust-region subproblem, so we can solve for $\lambda$ (and hence $w$) via root-finding with Newton's method. %described in~\cite{NocedalWright.optim.2000} that has a superlinear rate. 
The expensive step~\eqref{eq:tr} is implemented using Cholesky factors of the sparse matrix $\left(\frac{1}{\gamma} (L^TL) + \frac{1}{\eta} I + \lambda A^TA\right)$. 
This system only changes when $\eta$ is updated. 

\subsection{Alternative Approaches for Low-Rank \& Smooth Inversion}
\label{sec:alt}

There are alternative ways to model problem~\eqref{eq:MorozovSFR}. 
%All of them require a Tichonoff- or Morozov-type formulation. 
We consider two formulations and algorithms to compare with Algorithm~\eqref{alg:block-descent-smooth}.

\paragraph {Nuclear-norm formulation using FISTA. }
A simple convex formulation that uses smoothness,  data misfit, and a rank proxy (nuclear norm) is given by   
\begin{equation}
    \min_X \frac{\lambda}{2}\|\cA(X) - b\|^2 + \frac{1}{2\gamma}\|\cL(X)\|^2 + \|X\|_* \label{eq:fista}
\end{equation}
%(note that parameters are placed to keep notation consistent between formulations).
Formulation~\eqref{eq:fista} is the sum of a smooth and a simple function, and can be solved using projected gradient or Fast Iterative Shrinkage-Thresholding Algorithm (FISTA)~\cite{beck2009FISTA}, detailed in Algorithm~\ref{alg:fista}. This class of algorithms can be viewed as an extension of the classical gradient algorithm and is attractive due to its simplicity. The step size $\alpha$ is the reciprocal of the largest singular value of 
$(\lambda\cA^*\cA + \gamma^{-1} \cL^T\cL)$, and the operator $S_{\alpha}$ is the soft-thresholding operator:
\[
S_{\alpha}(\Sigma)_{ii} = \max(0, \Sigma_{ii} - \alpha). 
\]
\begin{algorithm}[H]
\caption{FISTA for~\eqref{eq:fista}.}
\label{alg:fista}
\begin{algorithmic}[1]
\State{\bfseries Input:} $X^0 = X^{-1} \in \R^{m\times n}, \ t^0 = t^{-1} = 1$
\State{Initialize: $k=0$}
\While{not converged}
\Let{$Y_{k}$}{$ X_k + \frac{t_{k-1}-1}{t_k}(X_k - X_{k-1})$}
\Let{$G_k$}{$Y_k - \alpha\lt((\lambda \cA^*\cA + \gamma^{-1} \cL^T\cL)\vecp(Y_k) - \lambda \cA^*b\rt)$}
\Let{$U, \Sigma, V^T$}{$\texttt{svd}(G_k)$}
\Let{$X_{k+1}$}{$US_{\alpha}(\Sigma)V^T$}
\Let{$t_{k+1}$}{$\frac{1+\sqrt{1+4(t^{k})^2}}{2}$}
\Let{$k$}{$k+1$}
\EndWhile
\State{\bfseries Output:} $X_k$
\end{algorithmic}
\end{algorithm}
Algorithm~\ref{alg:fista} uses gradients of the smooth terms, which requires appying $\cA$, $\cL$ and their adjoints, and the prox operator of $\|\cdot\|_*$, which requires thresholding on singular 
values computed via SVD (steps 6,7). These steps become prohibitively expensive as the dimensions of 
$X$ grow.

\paragraph{Smooth Factorized Formulation with L-BFGS.} 
To avoid the SVD steps of Algorithm~\ref{alg:fista}, we use the factorization strategy described 
in Section~\ref{sec:factform}:
\begin{equation}
\label{eq:lbfgs}
\min_{L,R}  \frac{\lambda}{2}\|\cA(LR^T) - b\|^2 + \frac{1}{2\gamma}\|\cL(LR^T)\|^2 + 
\frac{1}{2}\|L\|_F^2 + \frac{1}{2}\|R\|_F^2. 
\end{equation}
Formulation~\eqref{eq:lbfgs} is smooth with respect to the decision variables $L$ and $R$,
and at larger scales, the limited memory BFGS (L-BFGS)~\cite{NocedalWright.optim.2000} 
algorithm is a reasonable choice.

In the next section, we describe the travel-time interpolation problem for regional seismology, 
evaluate low rank and smooth regularization, and compare the performance of Algorithm~\ref{alg:block-descent-smooth} for~\eqref{eq:MorozovSFR}
to those of FISTA (Algorithm~\ref{alg:fista}) on~\eqref{eq:fista} and L-BFGS on~\eqref{eq:lbfgs}, in terms 
of computational efficiency and quality of reconstruction.

\section{Interpolation of synthetic travel time iMUSH data}
\label{sec:data}
\subsection{iMUSH Project and Data }
Mount St. Helens is the most active volcano in the Cascades arc. The Imaging Magma Under St. Helens (iMUSH) project was designed to illuminate the magmatic system beneath the volcano using a variety of geophysical and petrological methods, including active source, local earthquake, and ambient noise seismic tomography. %receiver functions, magnetotelluric imaging, and geochemical sampling and modeling. 
For the passive source seismic portion of the iMUSH project, 70 broadband seismometers were deployed from 2014 to 2016 within a 100km diameter circle around the mountain; these have an average station spacing of 10km, % [14 - what? ask], 
and are supplemented by permanent stations maintained by the Pacific Northwest Seismic Network (PNSN) and a temporary array of 20 broadband seismometers deployed by AltaRock Energy in June to November of 2016.% (IRIS network code YH, 2016).

During 2014-2016, over 400 earthquakes of magnitude 0.5 or greater occurred within 100km of Mount St. Helens. In addition, 23 borehole shots were set off during the summer of 2014 as part of the active source portion of iMUSH \cite{kiser2016vpvs}. These sources generated $P$- and $S$-waves which were recorded at the temporary and permanent networks by analysts, yielding over 12000 $P$-wave and 6000 $S$-wave arrival times. The arrival times were inverted to obtain 3D seismic velocity models and relocated hypocenters \cite{2017AGUFMcarl} %(right cite?) 
using the program struct3DP (written by Robert Crosson and based on previous codes used by Symons and Crosson, 1997).  % citation? 
The 3D $P$- and $S$-wave velocity models provide much greater detail than 1D velocity models used by the PNSN to locate earthquakes, or previous tomographic experiments in the area that focused on a smaller \cite{crosson1989helens,WAITE2009113} or broader scale \cite{moran1999rainier}.

The inversion system used by struct3DP is very large, with millions of parameters, and extremely scarce across the domain. It is regularized using an anisotropic 3D Laplacian operator which acts to smooth the velocity model in each spatial dimension, while allowing more variation in the vertical dimension. This tomographic inversion method depends on the overall quality of the data as well as the source-receiver raypath coverage. Data quality can vary depending on factors such as local noise sources (i.e. roads, rivers, animals, wind etc.), subsurface and near-surface path effects, distance from the hypocenter, attenuation, and source size. In the iMUSH local earthquake tomography, each travel time was assigned an uncertainty by an analyst based on characteristics of the seismic waveform. The raypath coverage depends on the location of both sources %which in the case of local earthquakes is out of human control (earthquakes occurred throughout the model area, but were concentrated to the west of MSH – include??), 
and seismometers, whose distribution is affected by road access in a rugged volcanic landscape, land ownership, winter snow depth, and station outages.

%from vandalism, equipment malfunction, and weather or animal issues.\\
Our goal here is to {\it increase} the raypath coverage for use in earthquake tomography by simulating seismometer locations and observations by interpolating noisy data. This is particularly useful for areas within the study region that did not have seismometers installed or which had seismometers out of service for extended periods of time. 

This experiment interpolates {\it synthetic} travel time residuals, which are calculated with respect to travel times predicted by a 1D velocity model. Travel times vary on the order of tens of seconds, while travel time residuals are generally between -1 and 1 seconds. We obtain synthetic travel times by using the forward modeling portion of the struct3DP code, which uses a finite difference 3D eikonal equation solver \cite{vidale19903dfd,hole19953dfdreflection}. These travel times are calculated for the best-fit 3D model from the iMUSH local earthquake tomography and compared to travel times through the PNSN S4 \cite{pavlis1983vstruct} 
%(or C3? – check) 
1D velocity model to obtain the synthetic residual. Similarly, for observed travel times, we subtract the travel time predicted through the 1D model from the observed travel time to obtain the observed residual. 
% These, however, are not used within the interpolation scheme.

Here, we define the experimental data. While raw iMUSH/PNSN data exists at station locations around the mountain (see blue dots in Figure \ref{fig:stations}), these iMUSH/PNSN stations are not gridded. Approximations or projections of this raw data onto the grid have yielded poor results for all of the aforementioned algorithms. To circumvent this issue for the time being, we instead 
conduct a forward problem with the 3D eikonal equation solver and use the best-fit 3D model to
%conduct an inverse problem with the 3D eikonal equation solver on the raw data to 
generate synthetic results for uniformly gridded stations in a square of 70km to 165km (with origin being at latitute 45.2, longitude -123.7). For the rest of this paper, we refer to the travel-time residuals between these synthetic travel times and those predicted by the 1D PNSN S4 model as the \textit{true} data, or $X_{true}$. To generate observations for the interpolation schemes, we subsample this data down to 15\% of total grid coverage over all sources by picking synthetic gridded stations near raw iMUSH/PNSN stations. This is a subset of synthetic stations meant to represent the distribution of the real-world stations that recorded the event. This subsampled data is then corrupted with noise. For each station, we generate a standard deviation parameter from the uniform distribution (0.03-0.15)(s); this captures each synthetic station's inherent uncertainty. 
The deviation range is based on raw iMUSH uncertainties recorded over the two year period. Then, for each station's data, we add zero mean random gaussian noise generated using that station's deviation parameter to create the \textit{observed} data, or $X_{obs}$. 
\subsection{Description of source tensor construction}
In order to apply the methods of Section~\ref{sec:Relax}, we have to specify a matricization of the data. 
The matricization we use is derived from the receiver grid, which is represented as a 95km-wide mesh with 5 kilometer spacing centered near Mount St. Helens. Again, the station grid is centered near the mountain and roughly coincides with the 70 stations deployed for 2 years.
%Other grid choices can be used in the algorithm, but best results were found with this mesh. 
Each entry in the matrix represents a point on this uniform receiver grid. 
Missing entries are designated by zeros, while observed entries are represented with the travel time residuals relative to the 1D model. 
Our experiments in this paper focus mainly on synthetic residuals of the nonlinear 3D modeled data relative to the 1D model, 
which provide a `ground truth' dataset we use to evaluate and compare interpolation techniques. 

%The omission of real data has two reasons: 1) this paper mainly focuses on enunciating the optimization method, and 2) real data occurs off the grid, and a suitable projection operator proved elusive. 

% Our experiments use two kinds of data: 1) residuals of 3D modeled data relative to the 1D model, and 2) residuals of observed data relative to the 1D model. Modeled data are on the grid by definition. 
% Real data are not on the grid, and we use a linear 3rd order Lagrange polynomial interpolation % citation? 
%  to project observations onto the grid.
% Our experiments use two kinds of data: 1) residuals of 3D modeled data relative to the 1D model, and 2) residuals of observed data relative to the 1D model. Modeled data are on the grid by definition. 
% Real data are not on the grid, and we use a linear 3rd order Lagrange polynomial interpolation % citation? 
%  to project observations onto the grid. \\
% using a linear combination of the nearest points in 2D (- get carl to send you operator if we want to use it.) 
%  \\
\begin{figure}
\centering     %%% not \center
\subfigure[Spatial locations of the sources (`*') around Mount St. Helens ($\Delta$). The box represents the edges of the uniform $20\times 20$ station grid. The grid panel (c) is taken from the source marked with a red `*'.]{\label{fig:sourceloc}\includegraphics[width=65mm]{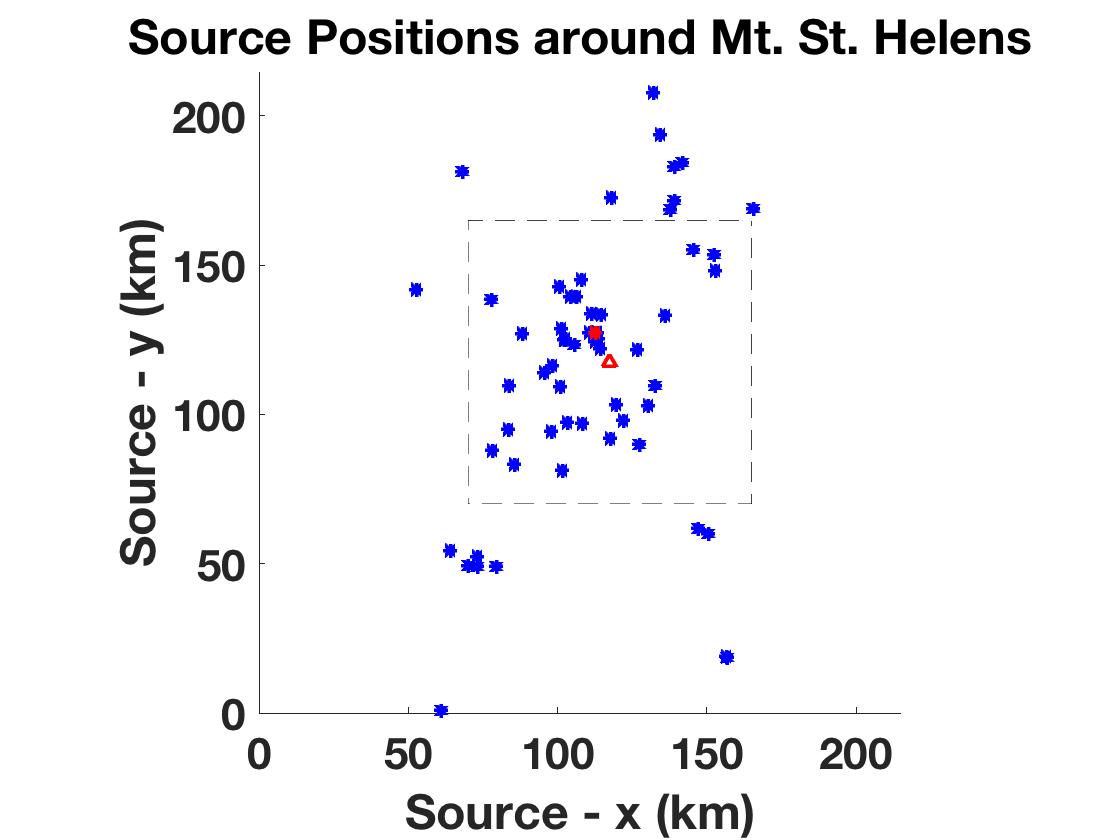}}\qquad
\subfigure[Receiver grid for a single source with the most datapoints. Black `*' represent synthetic stations on a grid, while blue `.' represent iMUSH/PNSN stations off-grid. ]{\label{fig:stations}\includegraphics[width=65mm]{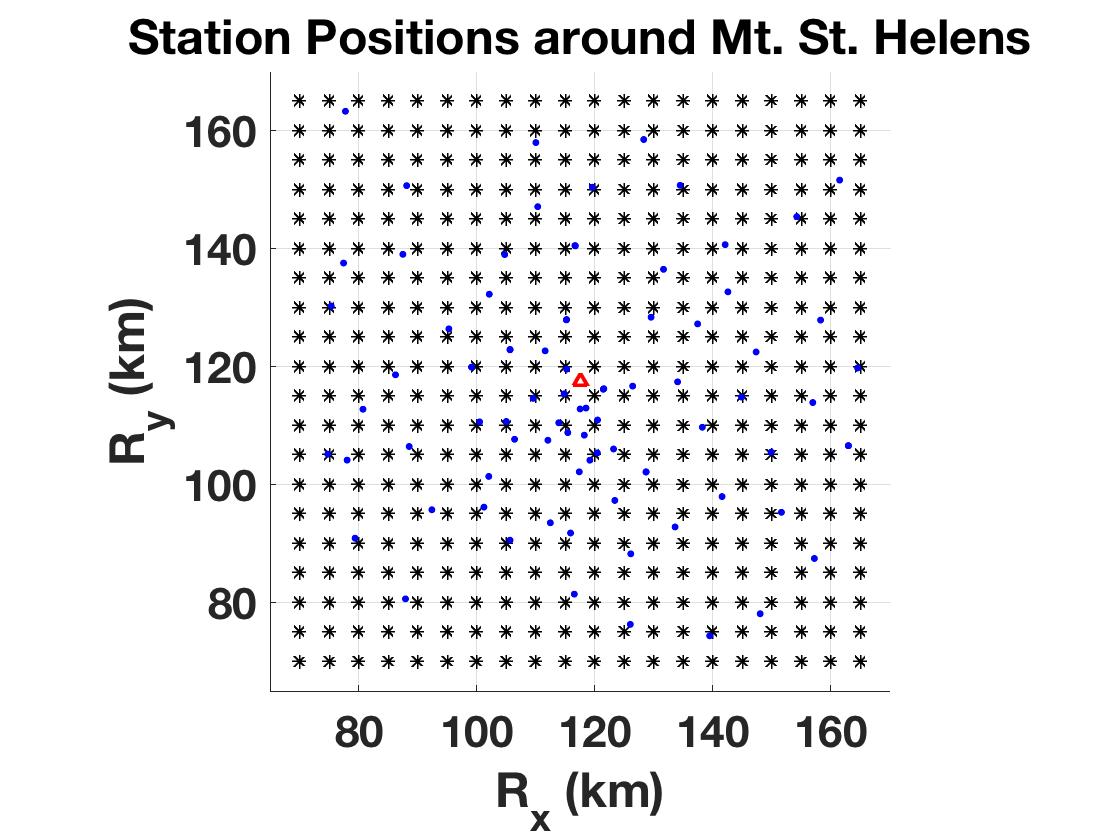}}
\subfigure[\textit{True} residual data for the receiver grid in panel (b). ]{\label{fig:rec_samp}\includegraphics[width=65mm]{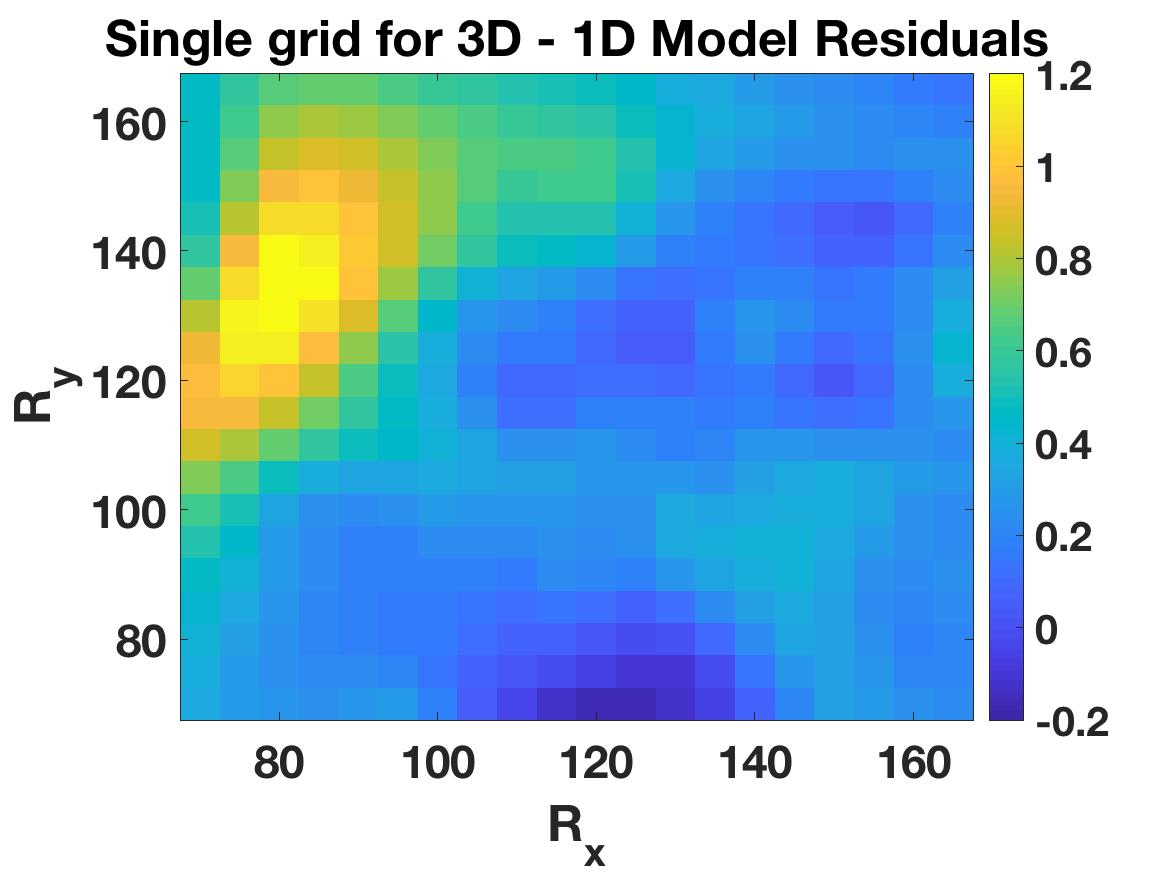}}
\caption{Data information for synthetic test.}
\label{fig:grid_info}
\end{figure}

%We use both global and local structure of the residual data. 
Each source represents a wave moving through the same media.
The underlying physical model suggests enforcing smoothness between station gridpoints. 
We further improve interpolation with low-rank methods by  
finding a tessellation of source-receiver grid matrices in which full data 
exhibits fast decay of singular values, while subsampled data does not. 
Intuitively, such a tessellation reflects the redundancy  of features across sources.  
Our combined goal is to interpolate $X$ from a subset of observations 
 by penalizing rank across sources, and nonsmooth local features.
 
 Each source $i$ at $(S_{x_i}S_{y_i})$ has an associated receiver grid of observations  
 $(R_x, R_y)\in (70, 165\text{km})\times (70, 165\text{km})$; each axis lies between 70km and 165km with 5km spacing. 
% This range was chosen for the concentration of sources around it (see Figure \ref{fig:sourceloc}), as well as the number of real stations located within it. %for which data is observed from that particular source. 
 The observation grid was chosen to lie close to the mountain and to contain a relatively large number of sensors  (see Figure \ref{fig:stations}). 
Initially, we observed 64 sources, where each source recorded at most 80 receivers (out of the potential 400). Source locations for the 64 sources are shown in Figure \ref{fig:sourceloc}, and a sample receiver grid for a particular source is given in Figure \ref{fig:stations} and Figure \ref{fig:rec_samp} is the \textit{true} data for that grid.  
The observed residual data is recorded in 3D tensor format with dimension $(n_{R_x}, n_{R_y}, n_s)$, where $n_{R_x}=20, n_{R_y} = 20, n_s = 64$ for the experiments.
Since observed residual data is a tensor, we have to choose a matricization (unfolding of the 3D tensor into a matrix) to exploit the induced low-rank structure. %is via matricizing the underlying tensor. 
%Matricization refers to a process of unfolding an n-D tensor into a 2-D matrix. Here, Low-rank structure refers to the small number of nonzero singular values or quickly decaying singular values of the underlying matrix. A challenge here is how to form the 2D matrix, $X$, from a 4D one. 
\begin{figure}
\centering     %%% not \center
\subfigure[Subsampled matricization $(R_x,R_y)\times(S_xS_y)$; the missing data (zeros) go across the rows.]{\label{fig:model_subsamp_form1}\includegraphics[width=65mm]{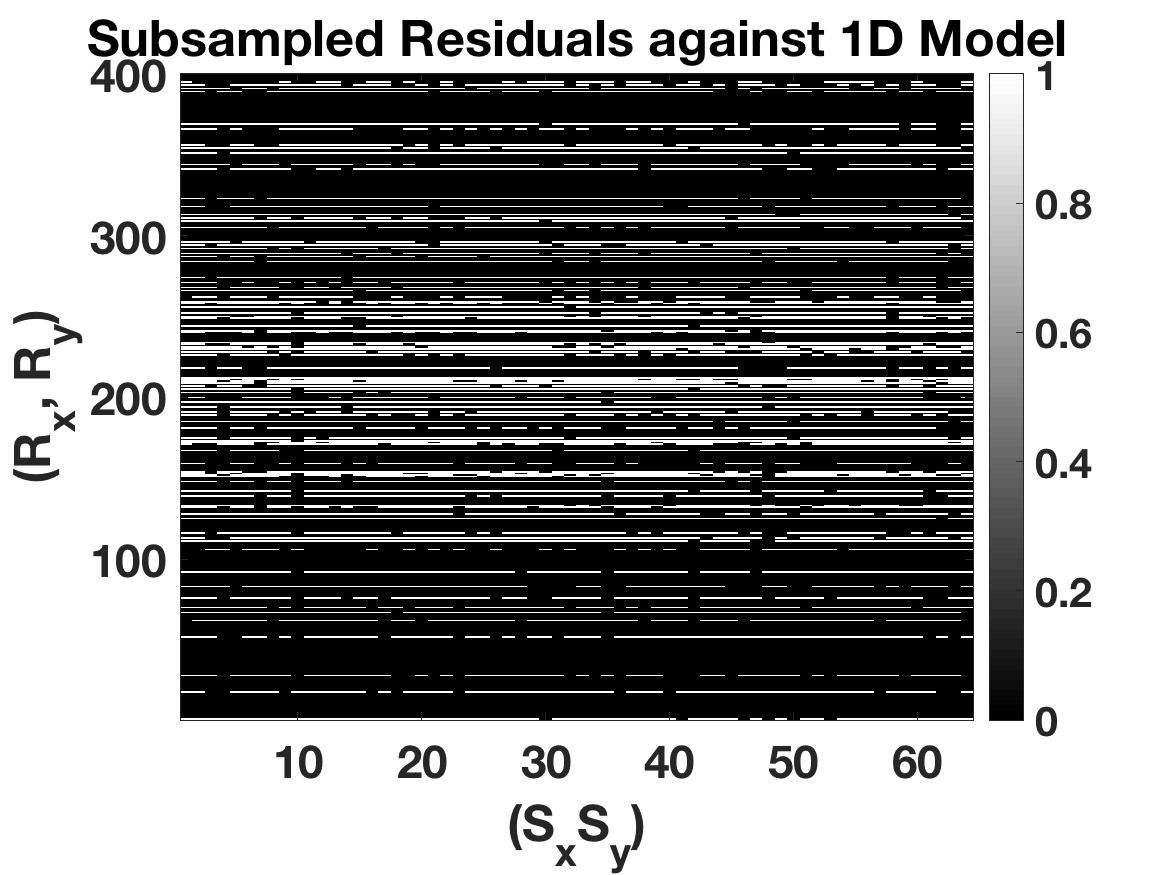}}\qquad
\subfigure[Subsampled matricization $(R_x, S_x)\times(R_y, S_y)$; missing entries (zeros) are are interwoven throughout the matrix.]{\label{fig:model_subsamp_form2}\includegraphics[width=65mm]{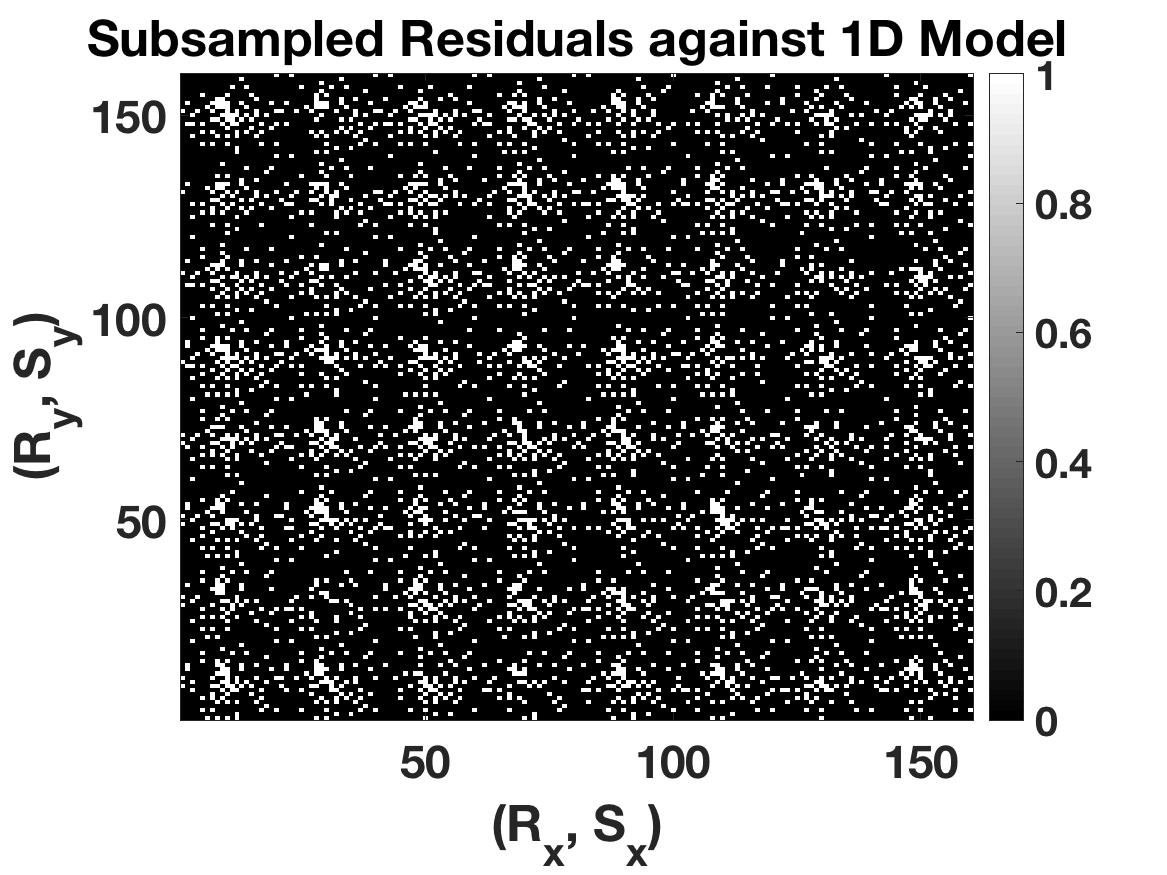}}
\caption{Different tensor formulations for low-rank interpolation.}
\label{fig:tensor_info}
\end{figure}

\begin{figure}[H]
\begin{center}
    \includegraphics[width=0.7\textwidth]{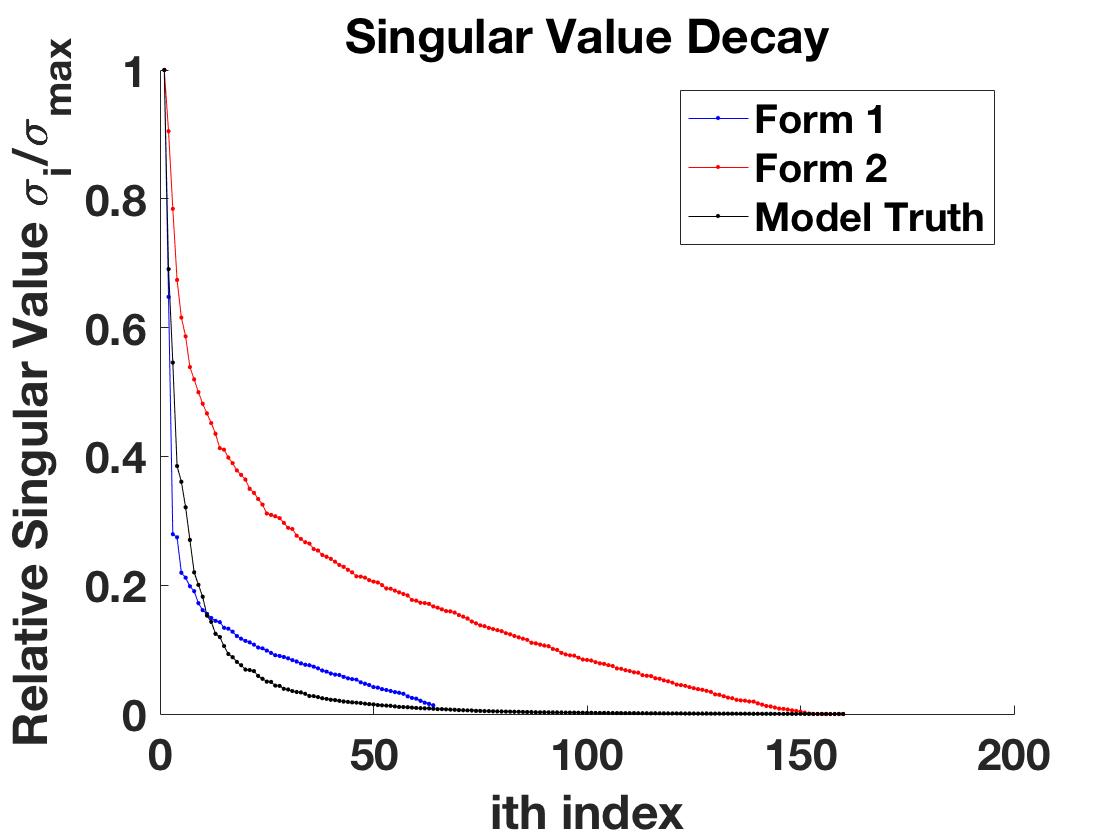}
    \caption{Singular value decay for matricization formulations 1 and 2 of the interpolation tensor. }
    \end{center}
    \label{fig:sval_decay}
\end{figure}

We consider  two such matricizations: % (as proposed in \cite{Curt} for seismic data): 
1) $(R_x,R_y)\times(S_xS_y)$, where we group receivers along columns and sources along rows, 
and 2) $(R_x,S_x)\times(R_y,S_y)$, where each receiver grid is block-inserted into the underlying matrix. On a 20x20 receiver grid with 64 sources, 
the first matricization $X\in\R^{400\times 64}$,  depicted in Figure \ref{fig:model_subsamp_form1},
%represented in tensor form by
%\begin{eqnarray*}
%            X = \begin{bmatrix} \Omega_{R_{x_1}R_{y_1}}^{S_{x_1} S_{y_1} } & \Omega_{R_{x_1}R_{y_1}}^{S_{x_2} S_{y_2} }&\dots & \Omega_{R_{x_1}R_{y_1}}^{S_{x_{64}} S_{y_{64}} } \\\\
%            \Omega_{R_{x_2}R_{y_2}}^{S_{x_1} S_{y_1} } & \Omega_{R_{x_2}R_{y_2}}^{S_{x_2} S_{y_2} }&\dots & \Omega_{R_{x_2}R_{y_2}}^{S_{x_{64}} S_{y_{64}} } \\
%            \vdots & \vdots & \dots & \dots \\
%            \Omega_{R_{x_{400}}R_{y_{400}}}^{S_{x_1} S_{y_1} } & \Omega_{R_{x_{400}}R_{y_{400}}}^{S_{x_2} S_{y_2} }&\dots & \Omega_{R_{x_{400}}R_{y_{400}}}^{S_{x_{64}} S_{y_{64}} } \\\end{bmatrix}
%\end{eqnarray*}
 %and 
 is obtained by letting $\Omega \in\R^{R_x,R_y} = \R^{400\times 1}$ be the vectorized receiver grid for single source $S_{x_i}S_{y_i}$ with $i = 1,\ldots, n_s$.
The second matricization $X\in\R^{160\times160}$, 
%written in tensor form by
%\begin{eqnarray*}
%	X = \begin{bmatrix} \Omega_{R_x\times R_y}^{ S_{x_1}S_{y_1} } & \Omega_{R_x\times R_y}^{ S_{x_2}S_{y_2} }&\ \dots & \Omega_{R_x\times R_y}^{S_{x_{8}}S_{y_{8}}} \\\\
%        \Omega_{R_x\times R_y}^{ S_{x_{9}}S_{y_{9}}} & \Omega_{R_x\times R_y}^{ S_{x_{10}}S_{y_{10}}} & \dots & \Omega_{R_x\times R_y}^{S_{x_{16}}S_{y_{16}}} \\
%        \vdots & \vdots & \dots & \dots \\
%        \Omega_{R_x\times R_y}^{ S_{x_{57}}S_{y_{57}}} & \Omega_{R_x\times R_y}^{ S_{x_{58}}S_{y_{58}}} & \dots & \Omega_{R_x\times R_y}^{ S_{x_{64}}S_{y_{64}}}\end{bmatrix}
%\end{eqnarray*}
depicted in Figure \ref{fig:model_subsamp_form2}, is obtained by letting $\Omega\in\R^{R_x\times R_y} = \R^{20\times 20}$ be the nested receiver grid for single source $S_{x_i}S_{y_i}$.
The subsampling scheme is crucial for the approach. The ideal situation is for the subsampled data to have high rank (slow decay of singular values), while the full data has low rank (fast decay of singular values). 
Then, it is feasible to recover the full volume by penalizing rank while matching observed data. 
 % then, the interpolation scheme will produce low-rank results. 
% Note that the first matricization is useful for subsampling, it is likely that an entire receiver position (row) is missing for every source. 
The first matricization does not satisfy this simple requirement: the subsampled matrix is itself low-rank and has rows and columns made up of zeros. Therefore we use the second matricization with tessellated sources, which indeed satisfies the requirement, see Figure \ref{fig:sval_decay}. 

In both matricizations, the sources are organized from high-energy to low-energy. In formulation 1, the sources are arranged such that the source with the greatest absolute residual values is in the first column and the least absolute residuals is in the last.
%\begin{minipage}\linewidth
%\centering
%\begin{minipage}[b]{.45\linewidth}

In formulation 2, the source with the greatest energy is in the top left, and the least in the bottom right; the energy then decays by source column (ie the next highest energy is the source immediately below in the row of sources). This produces an even slower decrease in singular values when receiver grid points are omitted. 

\section{Application to Synthetic Data}
\label{sec:results}
First, we evaluate the accuracy of using both smoothness and low-rank, compared to using either property alone. Then, we compare the speed and robustness of the new algorithm with competitors: 
 FISTA (Algorithm \ref{alg:fista}, Equation \ref{eq:fista}) and L-BFGS (Equation \ref{eq:lbfgs}). All tests use a tensored grid with the number of sources $n_s = 64$, and grid sizes $n_{R_x}, n_{R_y}=20$ for $R_x, R_y \in (70, 165)$ evenly spaced at 5km. The subsampling rate for every test is approximatestely $15$\%. We set the value $k$ (for the $LR^T$ formulation) to be 40 for all algorithms. 
 For all plots, north is up. 

First, we compare the `low-rank only' formulation~\eqref{eq:MorozovS} with the `smoothing only' formulation
\begin{align*}
	\min_W \|\cL(W)\|_2^2 \quad  \mbox{s.t.}\quad \|\cA(W) - b\|_2\leq \sigma.
\end{align*}
For this test, we use the inequality constraint $\|\cA(W) - b\|_2\leq \sigma$; the equality constraint ($\cA(W) = b$) is infeasible. % for smoothing only. %global (low-rank) and local (smoothing) only.
%The equality-constrained formulation above requires solving a single augmented linear system; using the constraint $\|\cA(W) - b\|_2 \leq \sigma$
This requires a root-finding algorithm nearly identical to the $W$-update of Algorithm~\ref{alg:block-descent-smooth} in Section~\ref{sec:Relax}.   
For the low-rank formulation, we use the matricization in Figure \ref{fig:model_subsamp_form2}. % to complete with low-rank interpolation. 
The convergence criteria ($l_2$-norm of minimized variable(s) iterate difference) for all algorithms is $10^{-10}$. All model  hyper-parameters are shown in Table \ref{tab:params}.
\begin{table}[H]
      \caption{\label{tab:params}Model hyper-parameters. We set $\gamma = 6.45\times10^{-7}$ for all algorithms except low-rank and smoothing only (where it is not used). 
    %It also incorporates the step-size taken by the tensored Laplacian matrix $\cL$. 
    The smoothing only formulation requires a root-finding problem for the inequality constraint, and has no Max Iteration value corresponding to block-coordinate-descent portion. 
    For FISTA, the step-size $\alpha = \|\cL\|_2^{-2}$, 
   the reciprocal of the largest singular value of $\cL$. In the joint formulation, we update $\eta$ every 30 iterations; in the low-rank only, we update every 100 iterations.}
  \begin{indented}
  \item[]\begin{tabular}{|l|l|l|l|l|l|l|l|l|l|l|l|l}
    \hline
    % \multicolumn{2}{c}{Error}                   \\
    % \cmidrule{2-3}
    Alg &  $\lambda$ & $\eta$ & $\eta_f$ & Max Iterations &$\sigma$ \\
    \hline
    Combined - VR Exact&   0.0111  & 0.5 &4.17 & 90 & 0 \\
    Combined - VR Noise&   0.0111  & 0.5 &4.17 &  90 & 3.719\\
	FISTA &   $2.2222\times 10^{-4}$   & --- & ---  & 1500 & ---\\
	L-BFGS &  $1.1111\times 10^{-4}$  &  --- &---  & 1500 & ---\\
	Smooth only&  --- & --- &---   & --- & 3.719  \\
	Low-rank only & --- &1.0  & 4.17 & 500 & 3.719\\
        \hline
  \end{tabular}
  \end{indented}
\end{table}
%although both algorithms essentially stop changing after approximately 50 iterations. 
For the proposed relaxation algorithm, 
$\eta$ is updated by $\eta_+ = \eta_f\eta$ for $\eta_f = \frac{\sum_{i}\hat\sigma_i}{k}$ at every 30th iteration. Here, $\hat\sigma_i$ are the singular values of the pre-interpolated matrix.
% \% = 9.38$. 
For low-rank only (with the $l_2$-norm), we update $\eta$ ever 100 iterations. Increasing $\eta$ can accelerate the algorithm~\cite{zheng2018fast}. 

\begin{figure}[H]
\centering     %%% not \center
\subfigure[Low-rank.]{\label{fig:lr_int}\includegraphics[width=60mm]{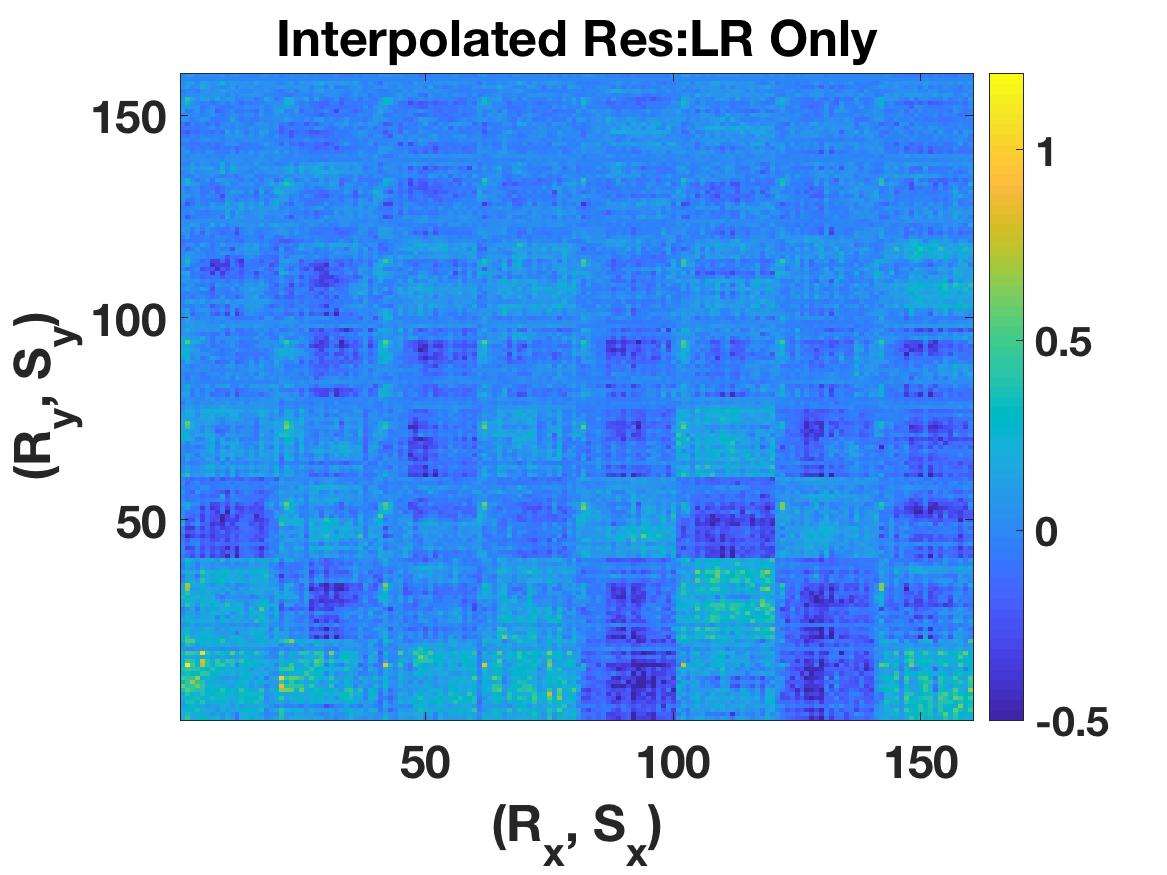}}
\subfigure[Smoothing. ]{\label{fig:sm_int}\includegraphics[width=60mm]{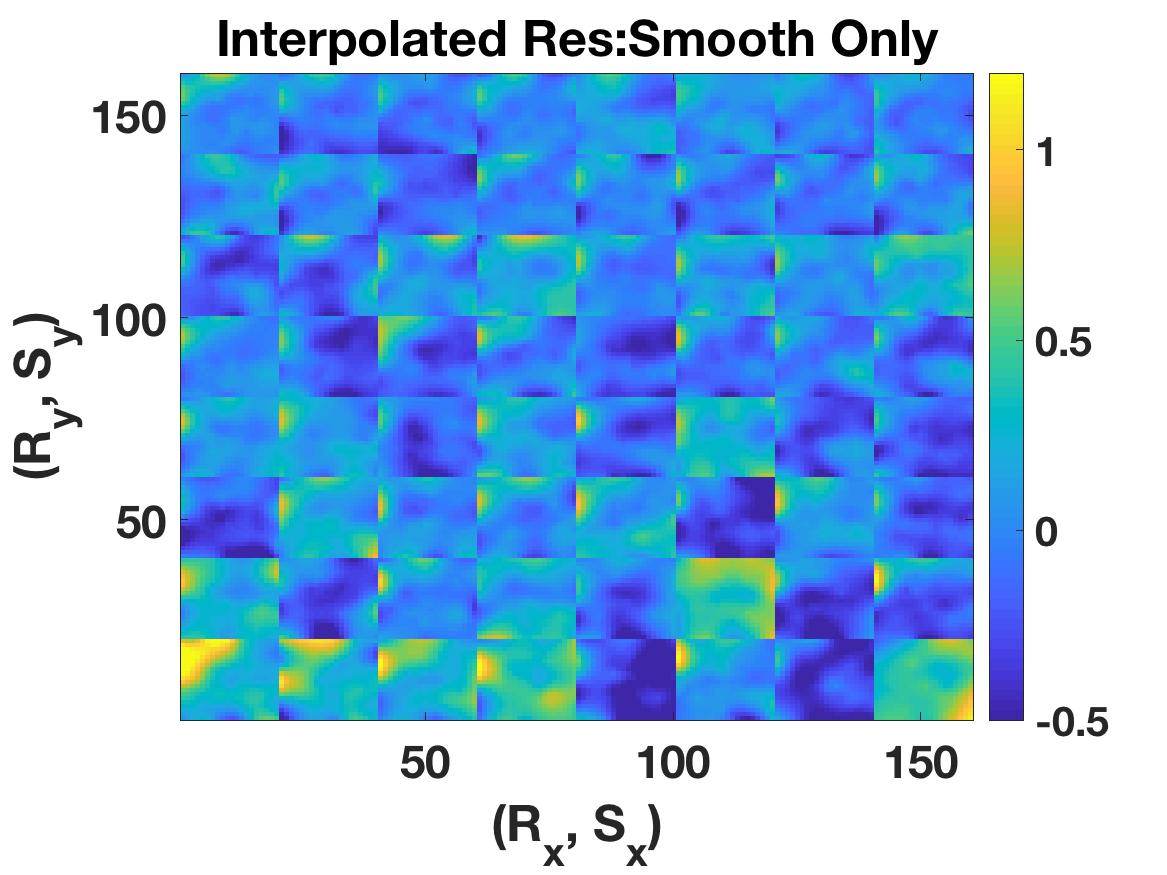}}
\subfigure[Observed model residuals, $X_{obs}$. ]{\label{fig:obs_all}\includegraphics[width=60mm]{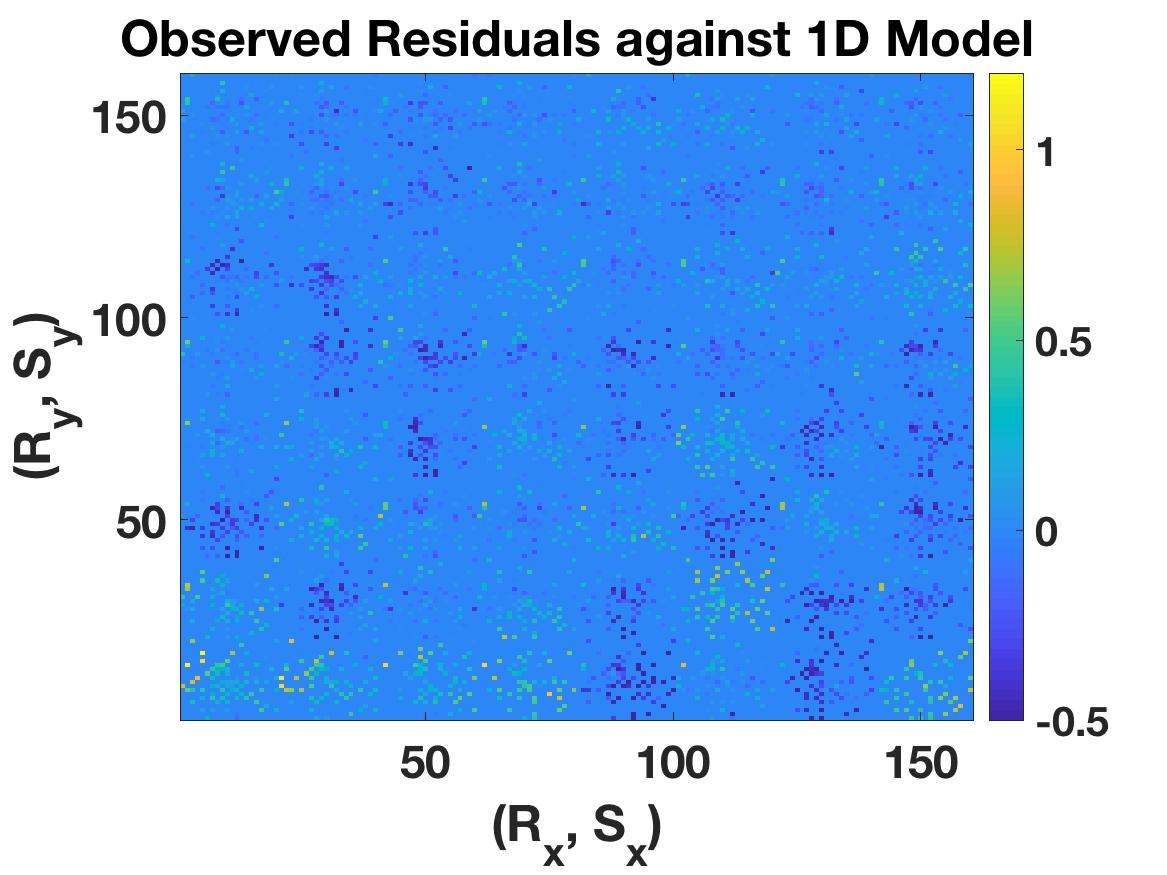}}
\subfigure[\textit{True} model residuals, $X_{true}$. ]{\label{fig:model_true}\includegraphics[width=60mm]{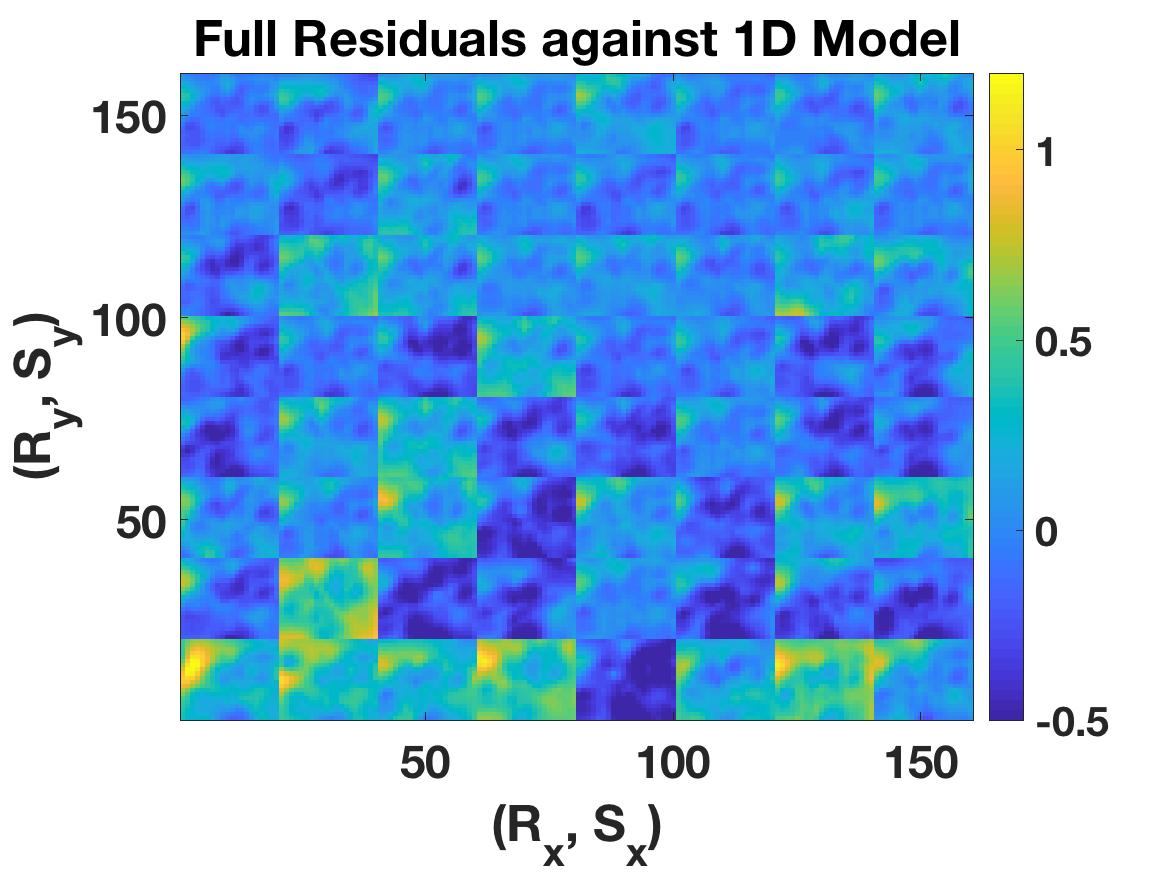}}
\caption{Full tensor results for Low-rank and smoothing only.}
\label{fig:lrsm_all}
\end{figure}

%Note that the value itself is somewhat arbitrary, but accelerating the convergence of $W$ and $LR^T$ produces better results here as well as in \cite{zheng2018fast} for certain problems.  
The interpolated results for these two schemes and the \textit{true} data are shown in Figure \ref{fig:lrsm_all} with root-mean-squared (RMS) errors for both observed and interpolated data listed in Table \ref{tab:results}. The RMS error is calculated with respect to \textit{true} data, and is split into two categories: RMS for \textit{true} data at locations that were observed (designated $\cA(X_{true})$ or obs), and RMS for \textit{true} data for locations that were not observed (given by the complement $\cA^\mathsf{c}(X_{true})$ or int). Recall that in this context, \textit{true} data refers to the residuals of the 3D model compared to the PNSN S4 1D velocity model uncorrupted by noise. The RMS is not weighted by uncertainties. 
While the low-rank only interpolation can accurately capture the observed data, it fails to reproduce the missing data. 

\begin{figure}[H]
\centering     %%% not \center
\subfigure[Low-rank.]{\label{fig:lr_int_single}\includegraphics[width=60mm]{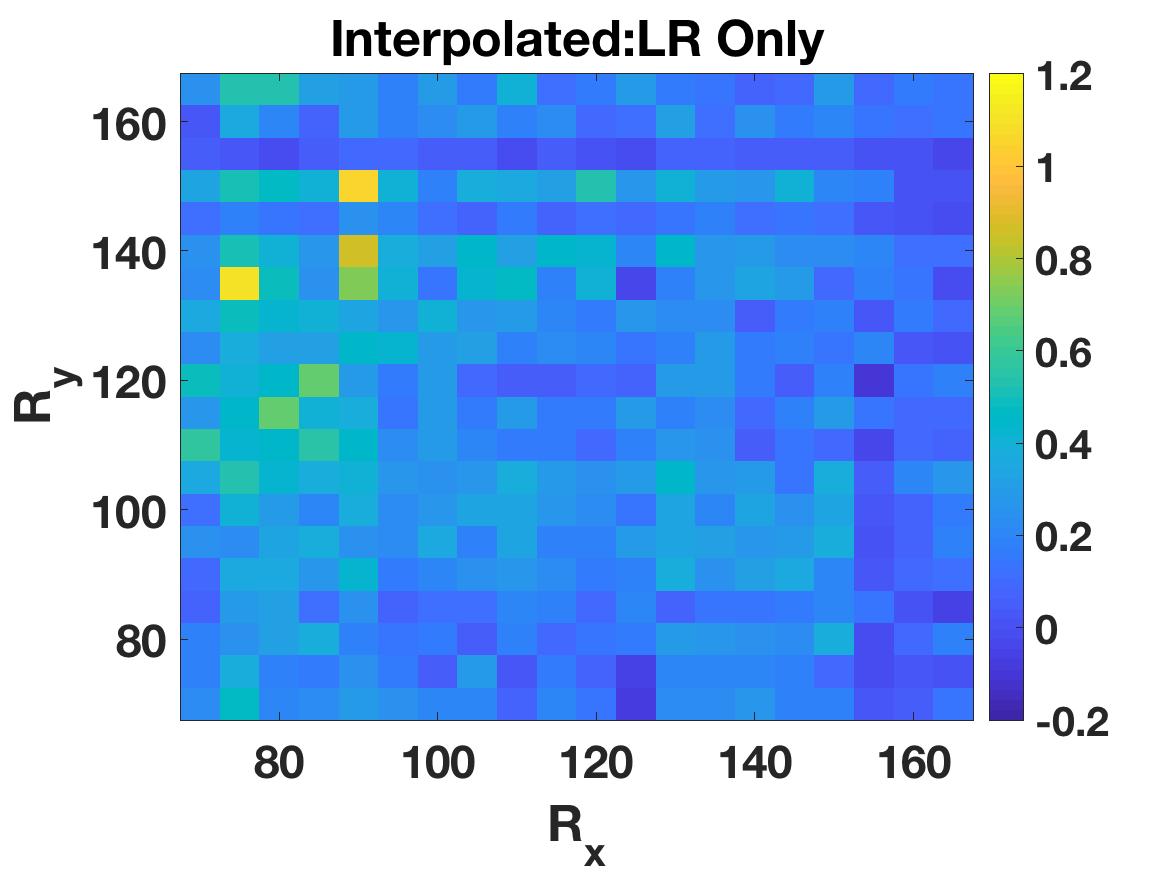}}
\subfigure[Smoothing only. ]{\label{fig:sm_int_single}\includegraphics[width=60mm]{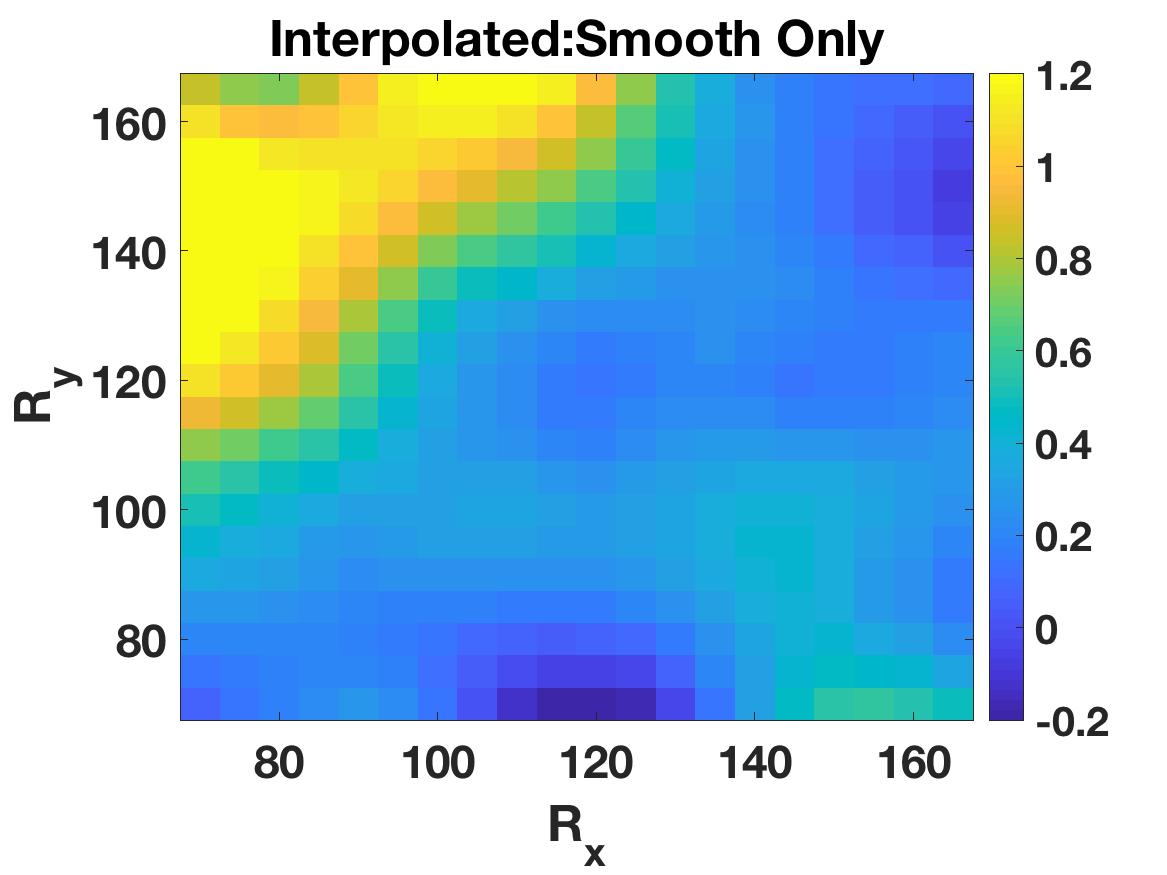}}
\subfigure[Observed model residuals. ]{\label{fig:obs_single}\includegraphics[width=60mm]{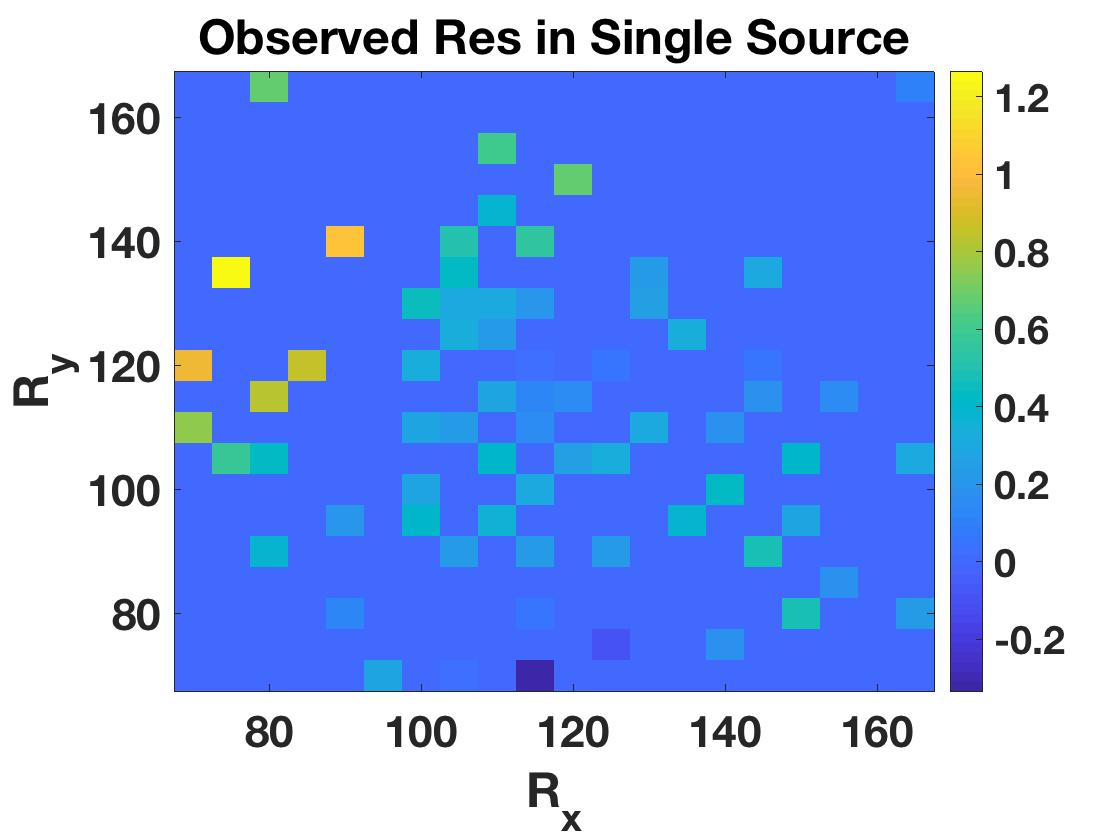}}
\subfigure[\textit{True} model residuals. ]{\label{fig:model_true_single}\includegraphics[width=60mm]{./figs/res_model_single1.jpg}}
\caption{Results for a single source with Low-rank and smoothing only.}
\label{fig:lrsm_single}
\end{figure}

\begin{figure}[H]
\centering     %%% not \center
\subfigure[Variable relaxation ($\sigma>0$).]{\label{fig:vr_diff_all}\includegraphics[width=60mm]{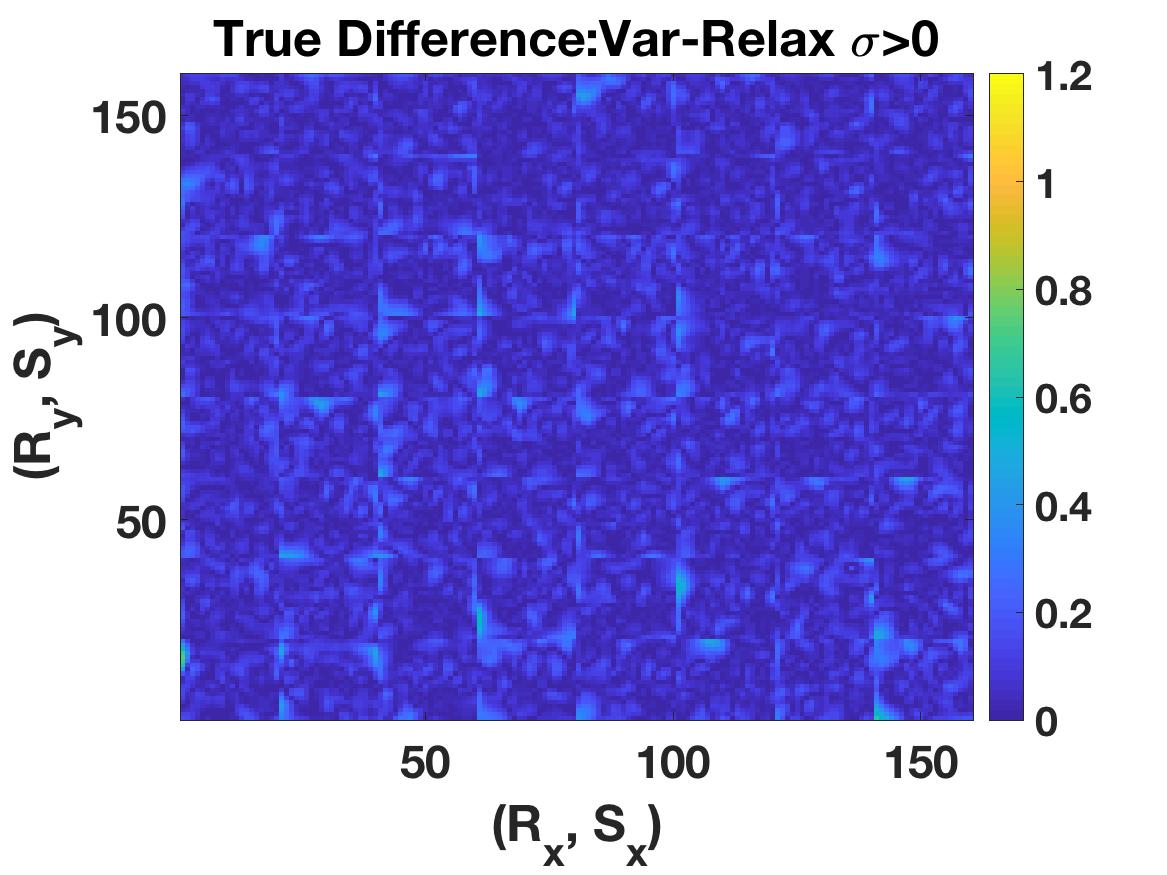}}
\subfigure[Variable relaxation ($\sigma=0$). ]{\label{fig:vrs0_diff_all}\includegraphics[width=60mm]{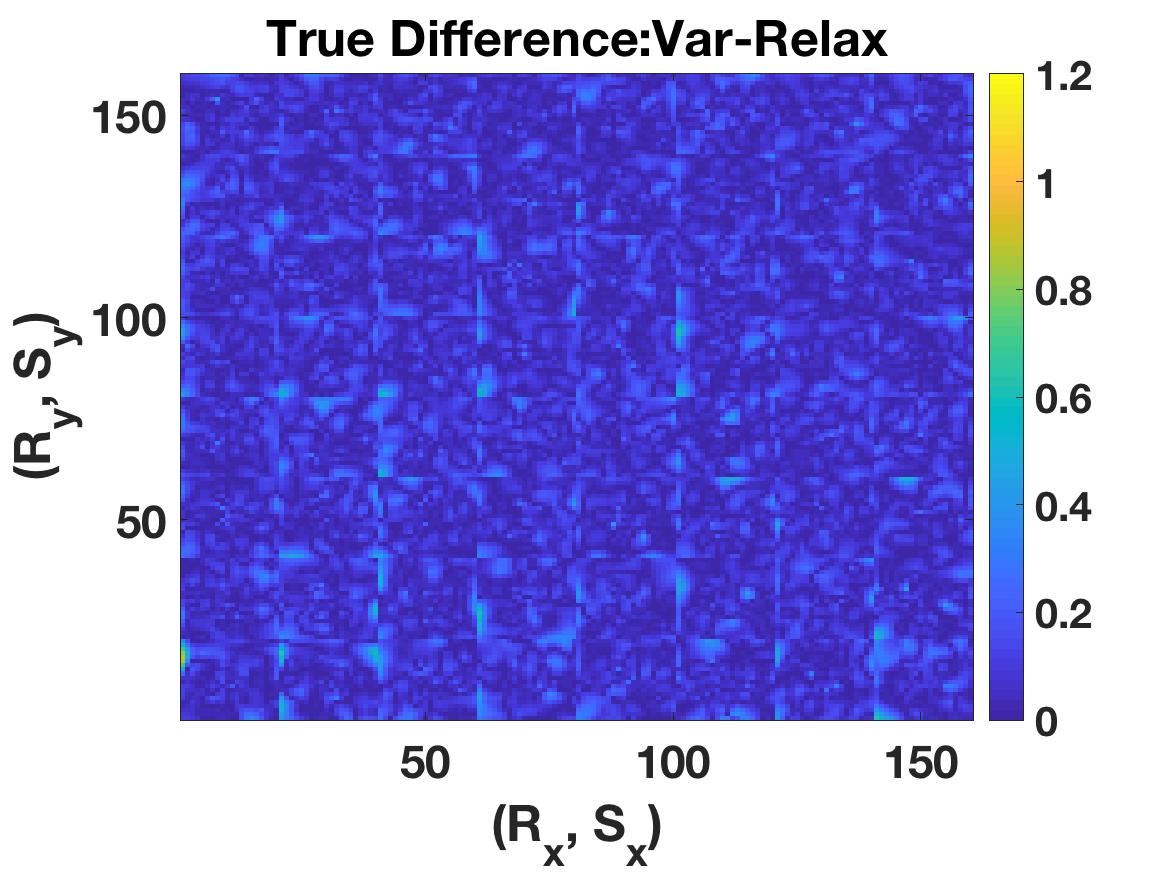}}
\subfigure[FISTA. ]{\label{fig:fista_diff_all}\includegraphics[width=60mm]{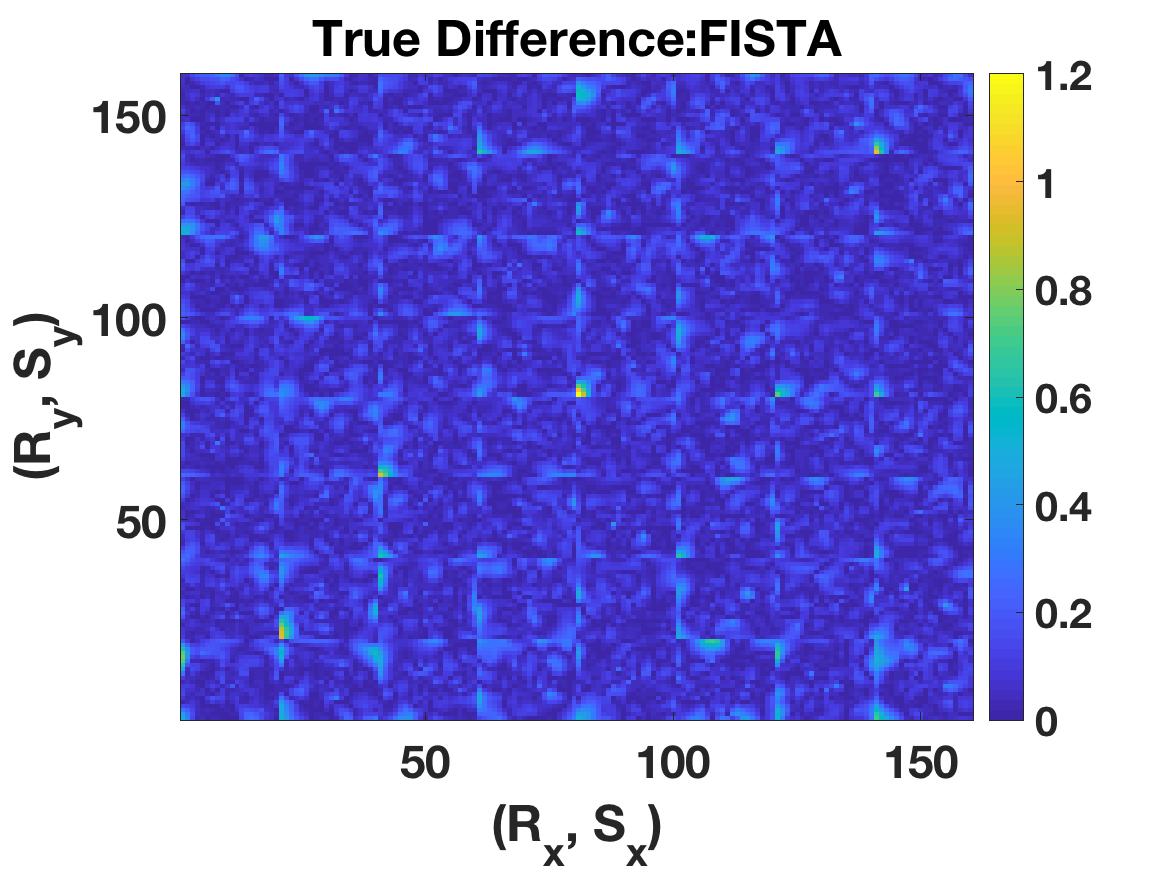}}
\subfigure[L-BFGS. ]{\label{fig:lbfgs_diff_all}\includegraphics[width=60mm]{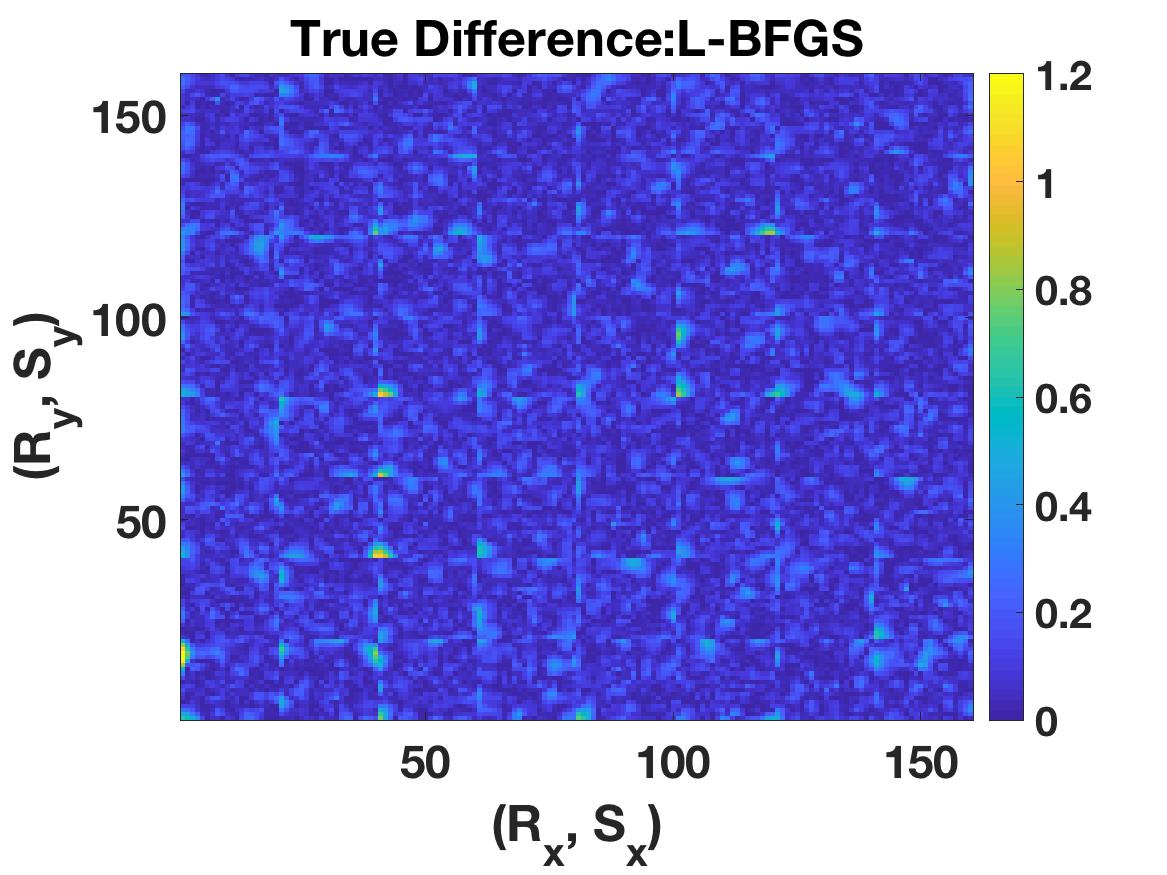}}
\subfigure[Low-Rank only.]{\label{fig:lr_diff_all}\includegraphics[width=60mm]{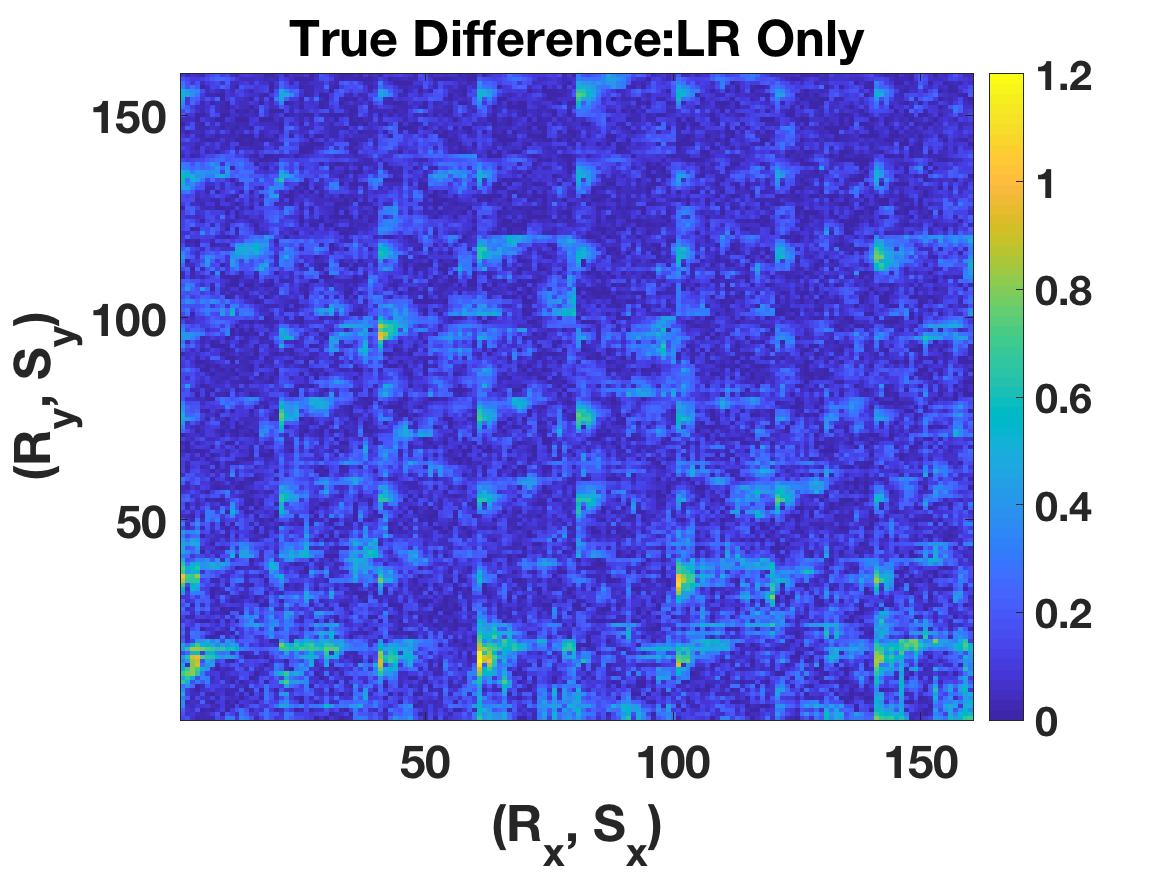}}
\subfigure[Smoothing only. ]{\label{fig:sm_diff_all}\includegraphics[width=60mm]{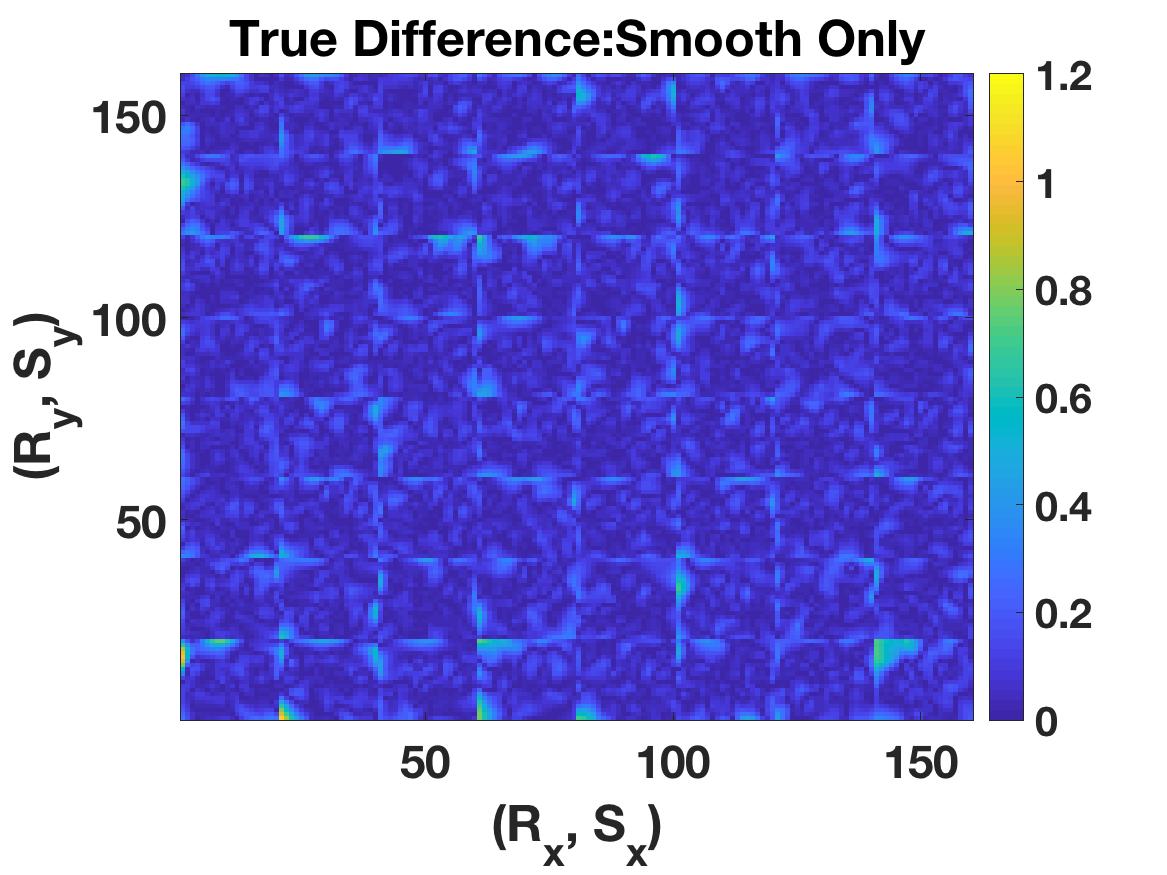}}
\caption{Full tensor $|X - X_{true}|$ for all algorithms. Note the significant scaling difference in the FISTA result.}
\label{fig:algcomp_diff_all}
\end{figure}

%Table \ref{tab:results} shows that the low-rank only produces a RMS that is comparable to the other algorithsm for observed entries, but fails entirely to capture the missing data. 
In Figure \ref{fig:lr_int}, the low-rank interpolation produces mono-color bands across the tensor: 
when entries for a receiver coordinate are missing in every source, the low-rank mechanism places zeros in the entire row or column. 
Likewise, if there is a single observed receiver in an entire row or column of the matrix, that value is propagated through every other entry in that row or column. 
The sparsity of the data impairs low-rank only interpolation,  which gives good results at sampling ratio of 70\% and above (these results not shown), but 
is less effective in our context. Smoothing alone, shown in Figure \ref{fig:sm_int} can capture major model dynamics, yet overestimates observed data worse than other methods and does not overall match data as well as other methods (see Table \ref{tab:results}). While the scale on Figure \ref{fig:sm_int} breaches is capped at 1.2, smoothing only actually has values larger than 1.5 in the bright yellow section, and the observed differences are also larger (Figure \ref{fig:sm_diff_all}). It also makes sense that smoothing alone would out-perform low-rank alone, since the underlying model used to generate data is inherently smooth. 
Some of the overall residual patterns are matched (larger residual energies are correctly placed around larger observed residuals), 
but the magnitude of the residuals is incorrect across all sources. 
A blown-up version of Figure \ref{fig:lrsm_all} for the source with the highest number of observed data ponts is shown in Figure \ref{fig:lrsm_single}. 
%Given the observed residuals in Figure \ref{fig:obs_single}, we can see that each formulation produces essentially for missing entries. 
%While low-rank alone has the capacity to capture observed entries, the goal of this procedure was to accurately predict sensor residuals in unknown locations. 
Both low-rank only and smoothing only formulations are inadequate for our application, especially with sparse data. \\
%The failure of each separate method to adequately capture the dynamics of the residuals demonstrates the utility of combining the local and global seismic information. \\


\begin{figure}[H]
\centering     %%% not \center
\subfigure[Variable relaxation ($\sigma>0$).]{\label{fig:vr_int}\includegraphics[width=60mm]{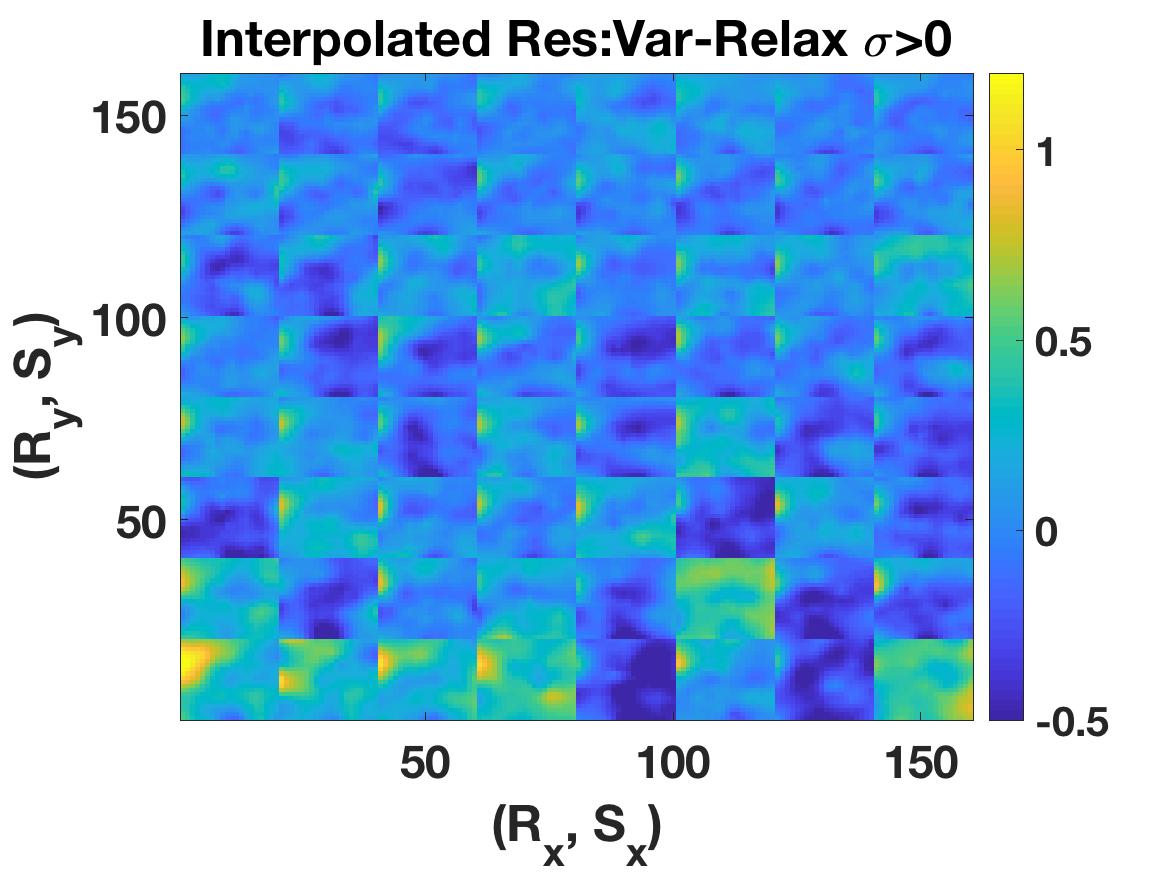}}\qquad
\subfigure[Variable relaxation ($\sigma=0$). ]{\label{fig:vrs0_int}\includegraphics[width=60mm]{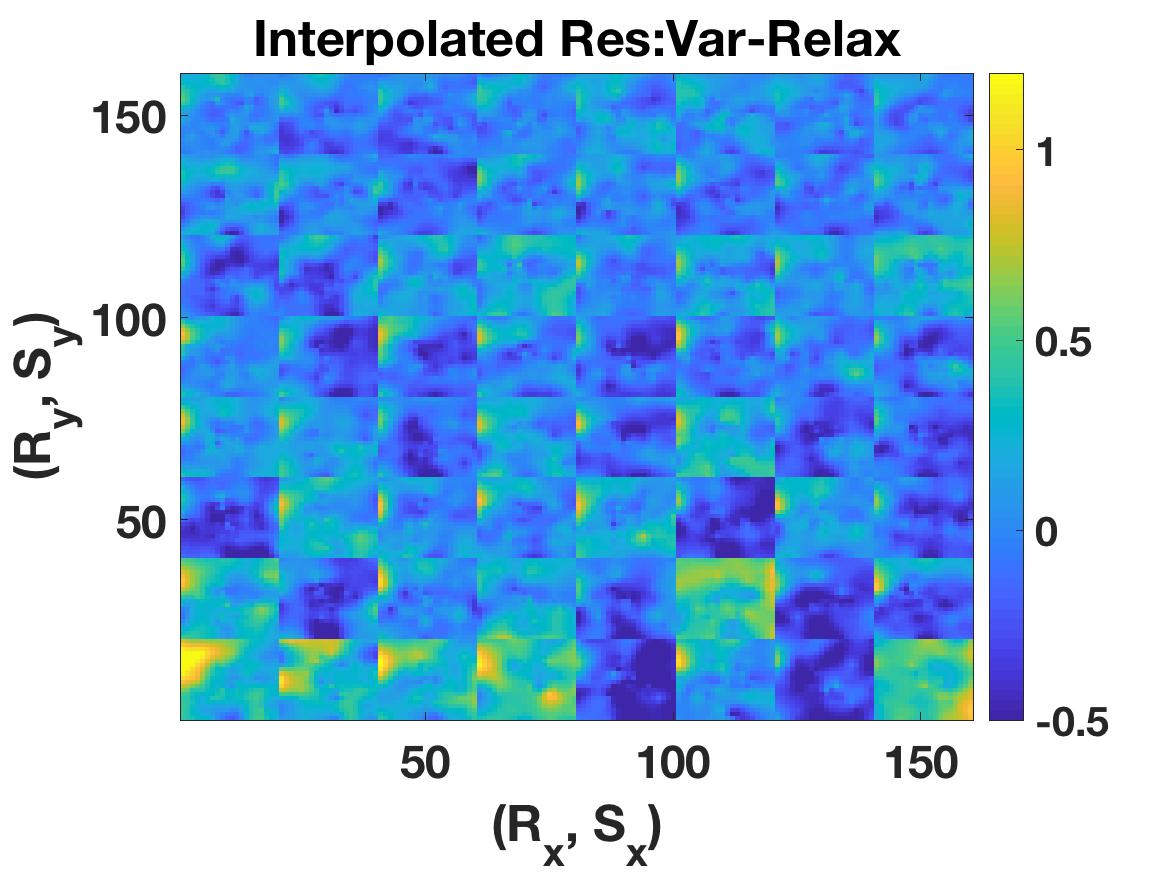}}
\subfigure[FISTA. ]{\label{fig:fista_int}\includegraphics[width=60mm]{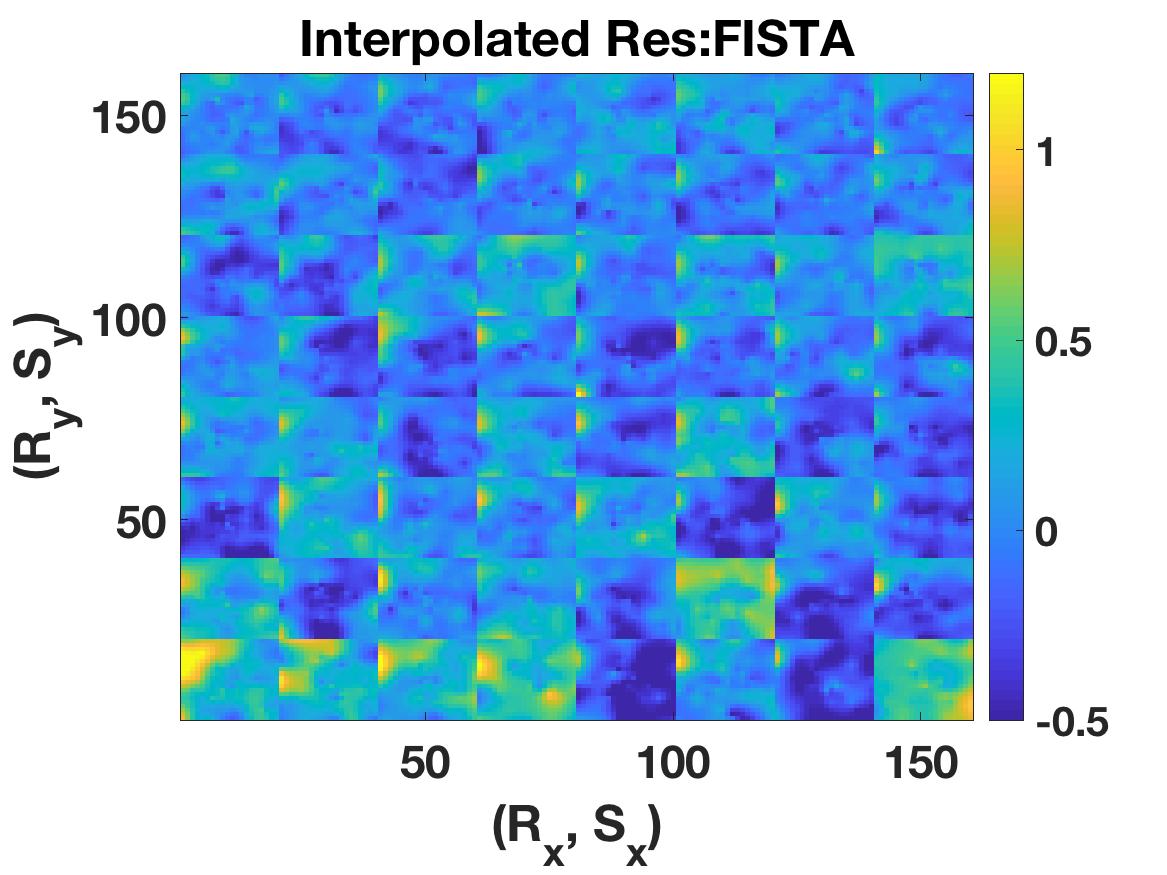}}\qquad
\subfigure[L-BFGS. ]{\label{fig:lbfgs_int}\includegraphics[width=60mm]{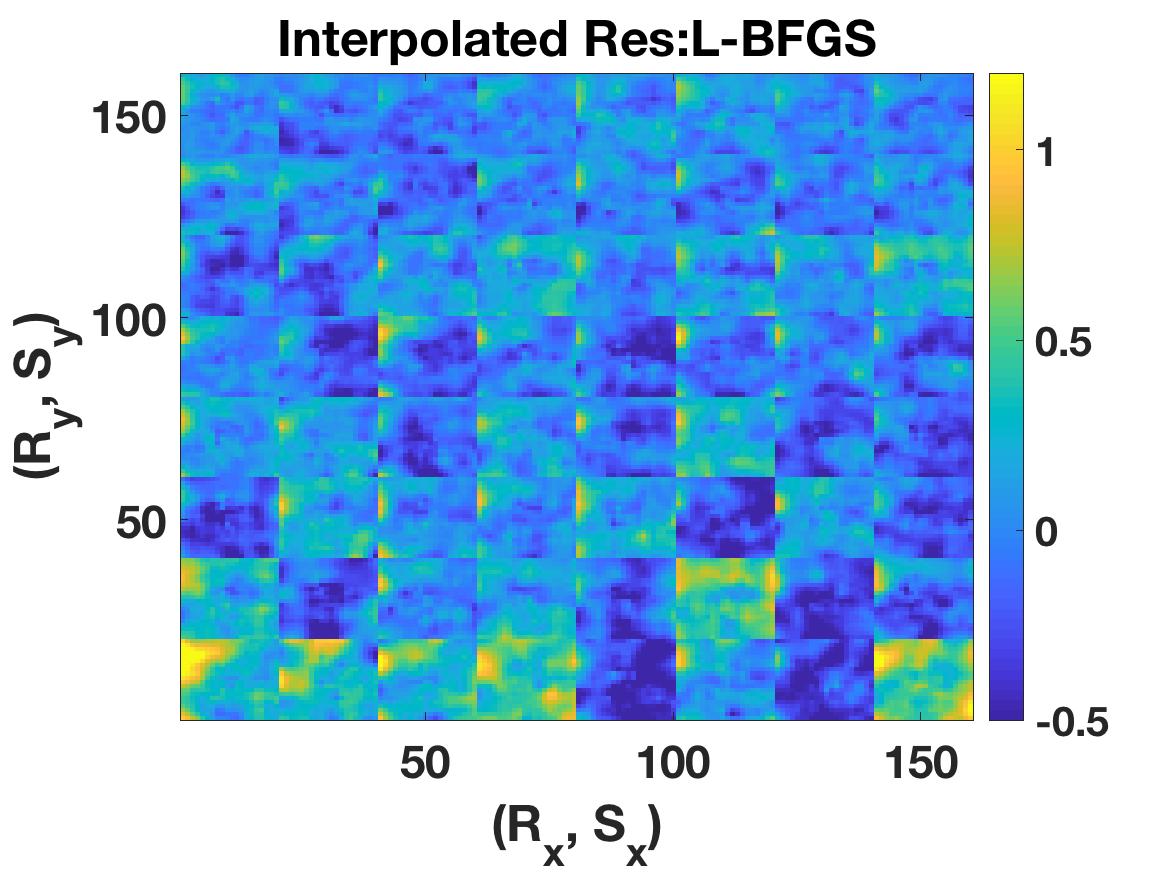}}
\subfigure[\textit{True} model residuals, $X_{true}$. ]{\label{fig:model_true_alg}\includegraphics[width=60mm]{./figs/res_model_all.jpg}}
\caption{Full tensor results for different algorithms and formulations.}
\label{fig:algcomp_all}
\end{figure}

\begin{figure}[H]
\centering     %%% not \center
\subfigure[Variable relaxation ($\sigma>0$).]{\label{fig:vr_int_single}\includegraphics[width=60mm]{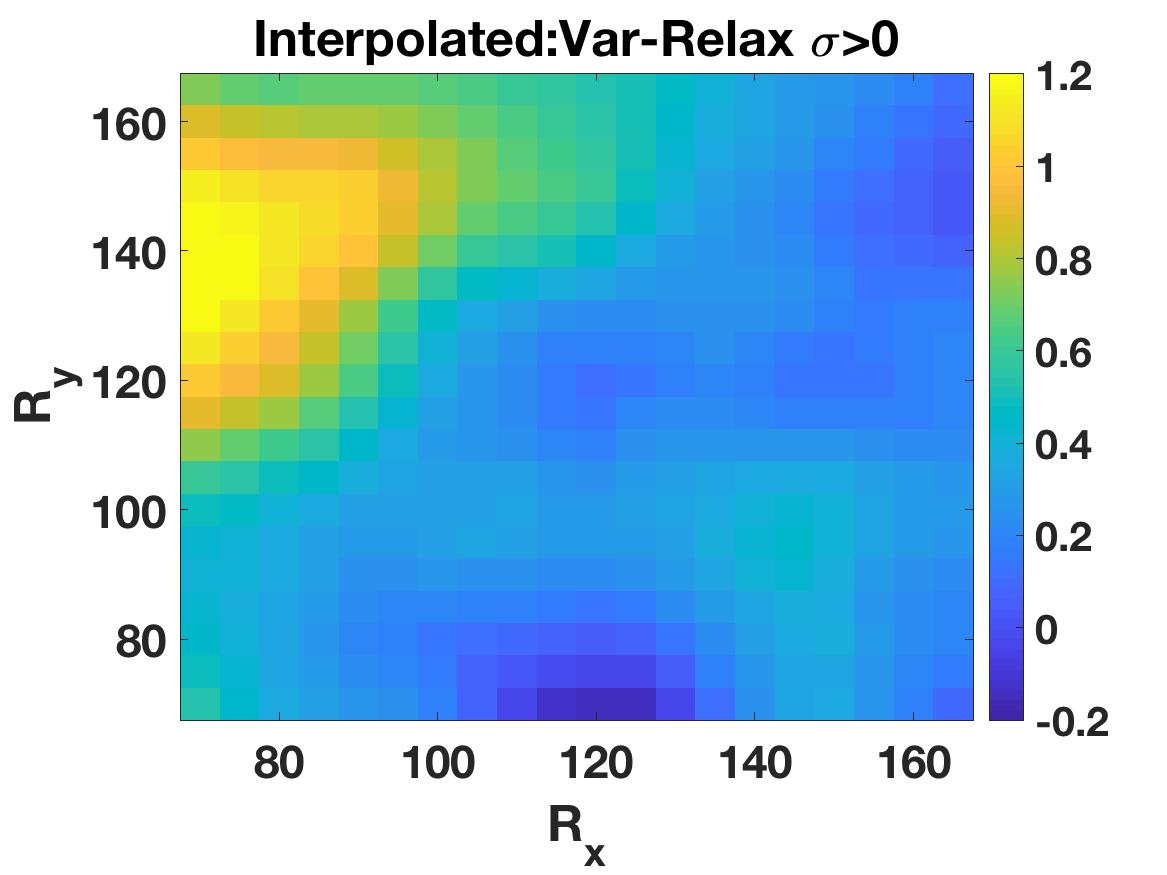}}
\subfigure[Variable relaxation ($\sigma=0$). ]{\label{fig:vrs0_int_single}\includegraphics[width=60mm]{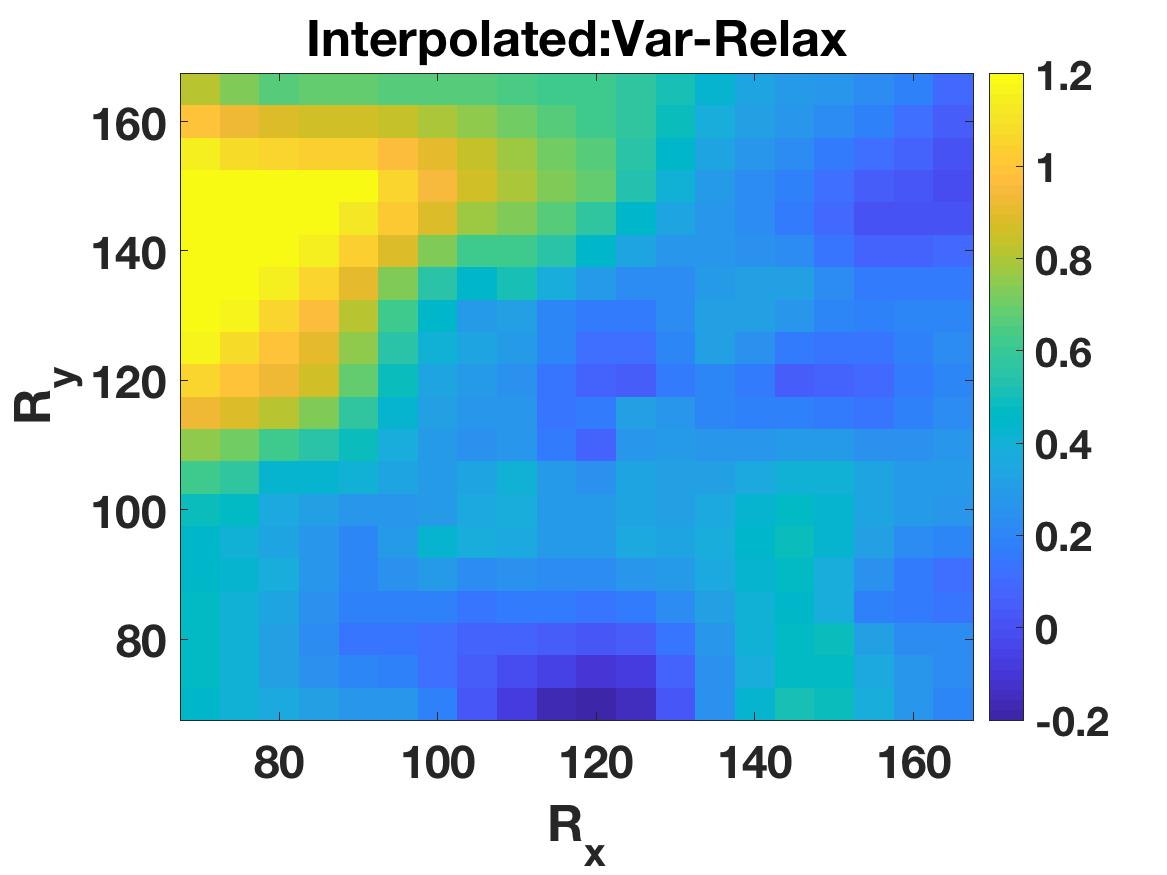}}
\subfigure[FISTA. ]{\label{fig:fista_int_single}\includegraphics[width=60mm]{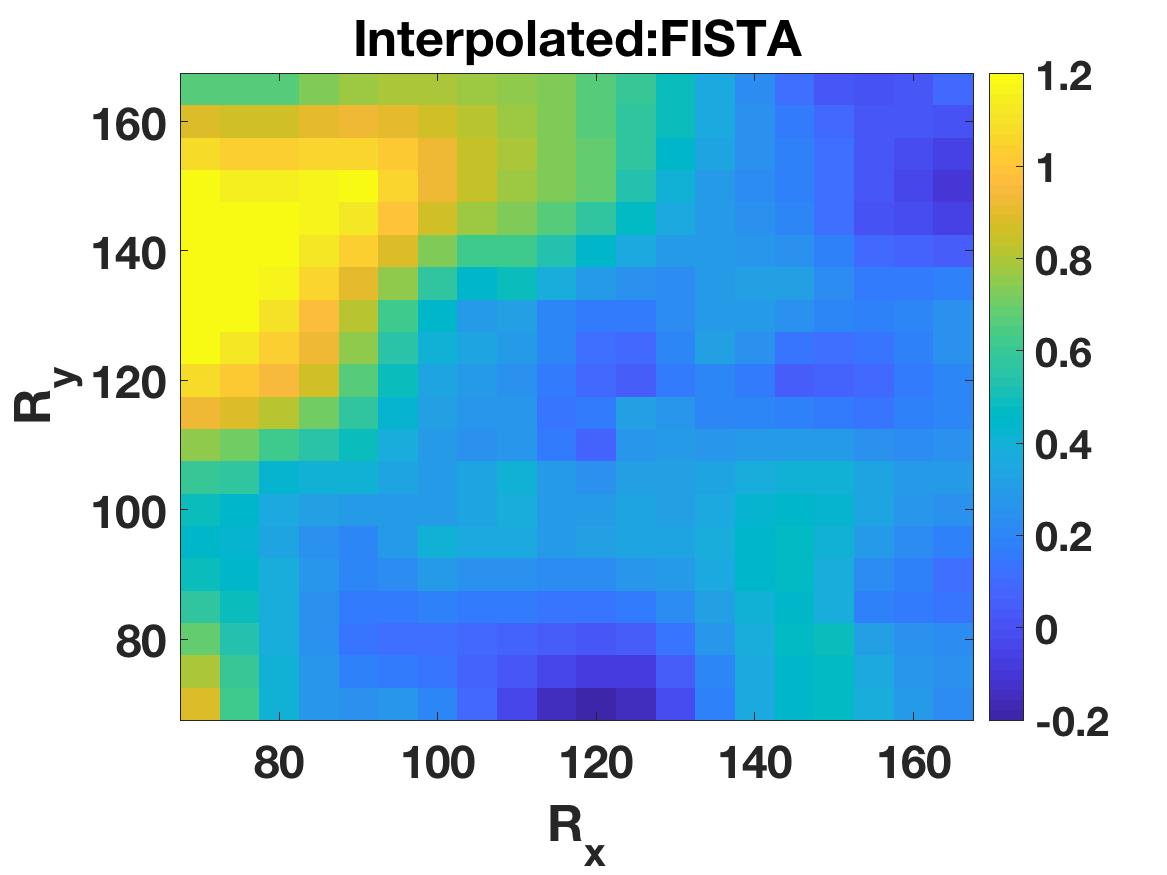}}
\subfigure[L-BFGS. ]{\label{fig:lbfgs_int_single}\includegraphics[width=60mm]{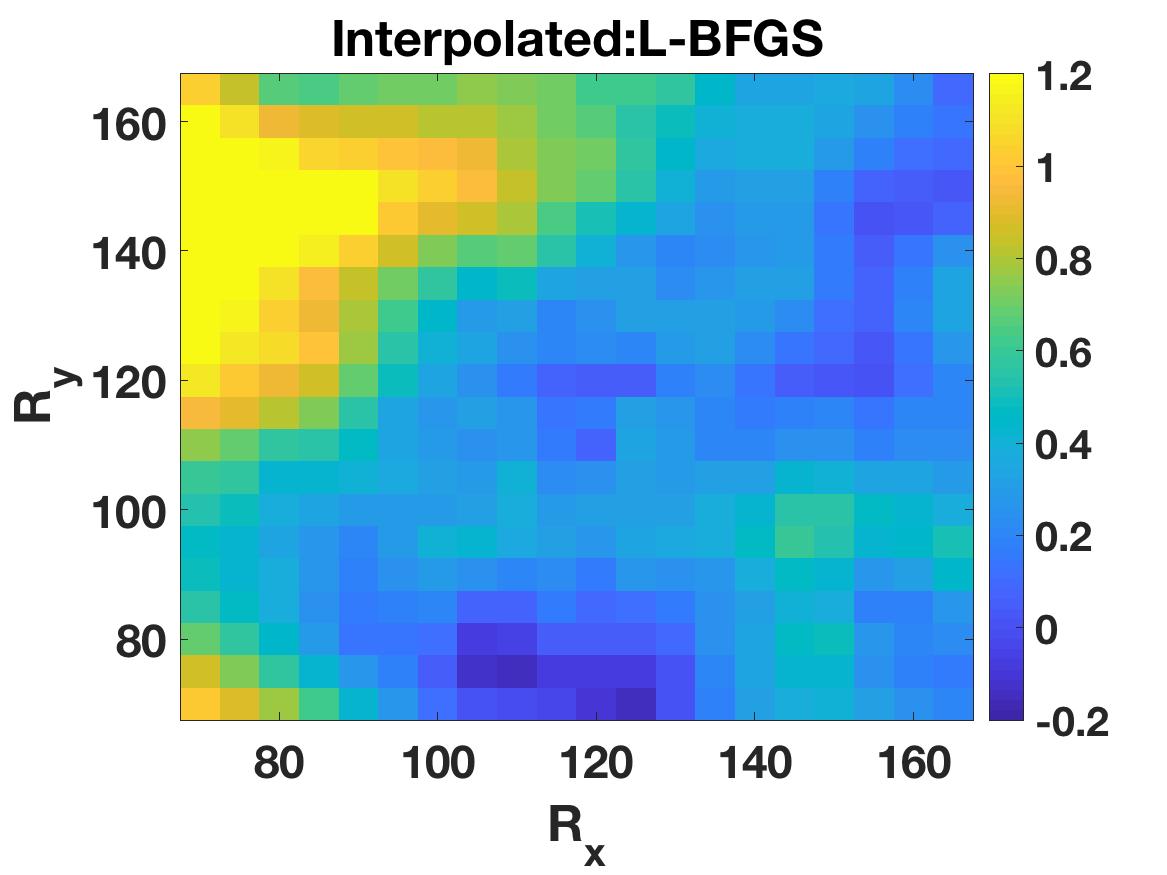}}
\subfigure[Observed model residuals. ]{\label{fig:obs_single2}\includegraphics[width=60mm]{./figs/obs_single.jpg}}
\subfigure[\textit{True} model residuals. ]{\label{fig:model_true_alg_single2}\includegraphics[width=60mm]{./figs/res_model_single1.jpg}}
\caption{Single source results for the different algorithms.}
\label{fig:algcomp_single}
\end{figure}

In contrast, we show in Figure ~\ref{fig:algcomp_all} that using both types local and global information can meaningfully interpolate the data. This can be done with
the new joint formulation as well as unconstrained formulations that use FISTA and L-BFGS algorithms; so we also test the efficacy of the new approach against these two competitors. 
%competing algorithms for interpolating noisy model data with the new algorithm. 
We consider two different choices of our variable relaxation algorithm: $\sigma = 0$ and $\sigma>0$. 
In the latter case, we assume that we do not know the true data misfit, and use available uncertainties to set $\sigma = \sqrt{\sum_{i = 1}^n0.06^2}\approx 3.72 (s)$ for $n$ being the size of $b$.
% using available uncertainties to set the data-fit level.  
%which uses the terminal feasibility of the low-rank only, and serves as a stand-in for the allowable misfit per observed datapoint (does this make sense?). 
The true data misfit is actually $\|b -\cA(X_{true})\|_2 = \sigma_{true} = 5.85(s)$ where $\cA(X_{obs})=b$ is our {observed} data. 
The results for tensored and single-source  matrix-completion are shown in Figures \ref{fig:algcomp_all} and \ref{fig:algcomp_single}, and in Table~\ref{tab:results}. 
The new approach (in bold) achieves better results than competing methods in similar amounts of compute time. 
%while also not relaxing or changing the defined problem in any way. 
Setting $\sigma>0$ in the variable relaxation scheme produces smaller RMS values, 
and in particular recovers missing data with higher accuracy. 
% 0.09 RMS error for FISTA/L-BFGS's observed values compared to VR's 0.097 RMS error for interpolated values. 
%In particular, the RMS for the relaxed approach on observed values is 0.06; so 
%we see better interpolation results are achieved by not fitting the data exactly. 
%when the observed  
% so we see that by \textit{relaxing} the data-misfit constraint purposefully, one can actually achieve better results than trying to fit the data exactly.\\
\begin{table}[H]
 \caption{\label{tab:results}Different formulations for model residuals with sampling rate of 15\%. Terminal feasibility is $\|\cA(X) - b\|_2 -\sigma$ at the algorithm's termination (which would mean $l_2$-norm data misfit where $\sigma=0$). Note that the feasibility is calculated against observed data while RMS is calculated against \textit{true} data. Recall that (obs) signifies $\cA(X_{true})$, while (int) stands for $\cA^\mathsf{c}(X_{true})$.}
  \begin{indented}
  \item[]\begin{tabular}{|l|l|l|l|l|l|l|l|l|l|l|l|l}
    \hline
    % \multicolumn{2}{c}{Error}                   \\
    % \cmidrule{2-3}
    Alg & Terminal Feasibility &Time (s) & RMS (obs) & RMS (int) \\
    \hline
    \textbf{Combined - VR Exact}& \textbf{0.0031} & \textbf{10.81} &   \textbf{0.09 }& \textbf{0.110} \\
    \textbf{Combined - VR Noise}& \textbf{1.18e-08}& \textbf{19.84}&  \textbf{0.06 }& \textbf{0.100} \\
	FISTA & 0.081 & 12.30  & 0.09 & 0.119\\
	L-BFGS &0.074 & 57.22 &0.09 & 0.128\\
	Smooth only& 0.0016 & 5.42  & 0.06 & 0.125 \\
	Low-rank only& 0.058 & 1.27  & 0.08 & 0.216\\
        \hline
  \end{tabular}
  \end{indented}
\end{table}
Table~\ref{tab:results} also shows the degree to which algorithms match the feasibility constraint. FISTA and L-BFGS use penalties rather than constraints, 
so the precise data-misfit level is hard to control. 
%is satisfied to some degree in each formulation. and solving a constrained problem with one of these methods would require one to relax the constraint by setting $\sigma=0$ and augmenting the cost function to take this data misfit into account (is this right?). 
The terminal feasibility is 0.08(s) for FISTA and 0.075(s) for L-BFGS, which means for each individual point, the values are close to fitting the observed entries. FISTA's feasibility level settles at approximately 0.08(s) and does not change after about 1000 iterations; similar behavior is seen in L-BFGS. 
%presumably as the low-rank and smoothing components of the cost function are continuing to decay. 
%A similar occurence happens for L-BFGS. 
Variable relaxation schemes can match the feasibility constraint to a high accuracy, 
with the explicit inequality constraint matching close to numerical precision. %\texttt{Note to sasha: our root-finding alg fails for $\sigma>5$.} \\

\begin{figure}[H]
\centering     %%% not \center
\subfigure[Variable relaxation ($\sigma>0$).]{\label{fig:vr_diff_single}\includegraphics[width=60mm]{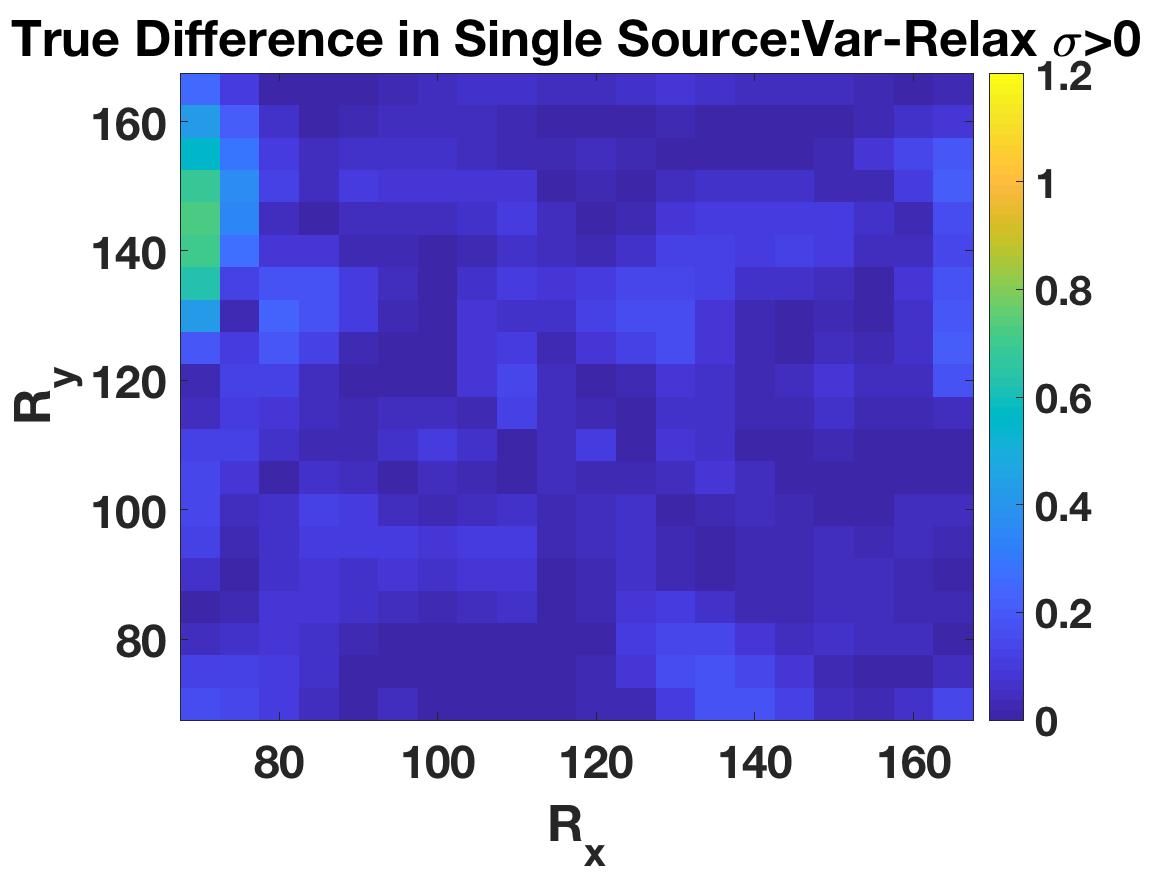}}
\subfigure[Variable relaxation ($\sigma=0$). ]{\label{fig:vrs0_diff_single}\includegraphics[width=60mm]{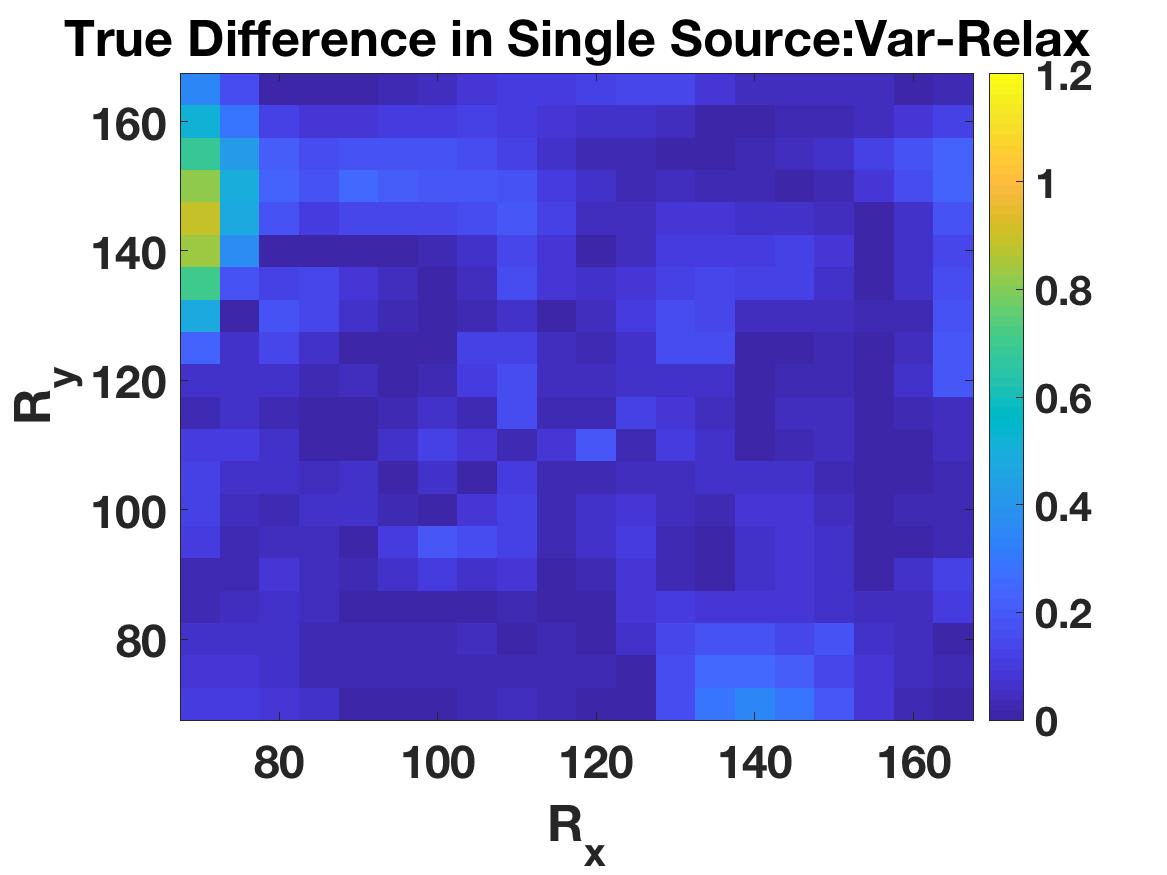}}
\subfigure[FISTA. ]{\label{fig:fista_diff_single}\includegraphics[width=60mm]{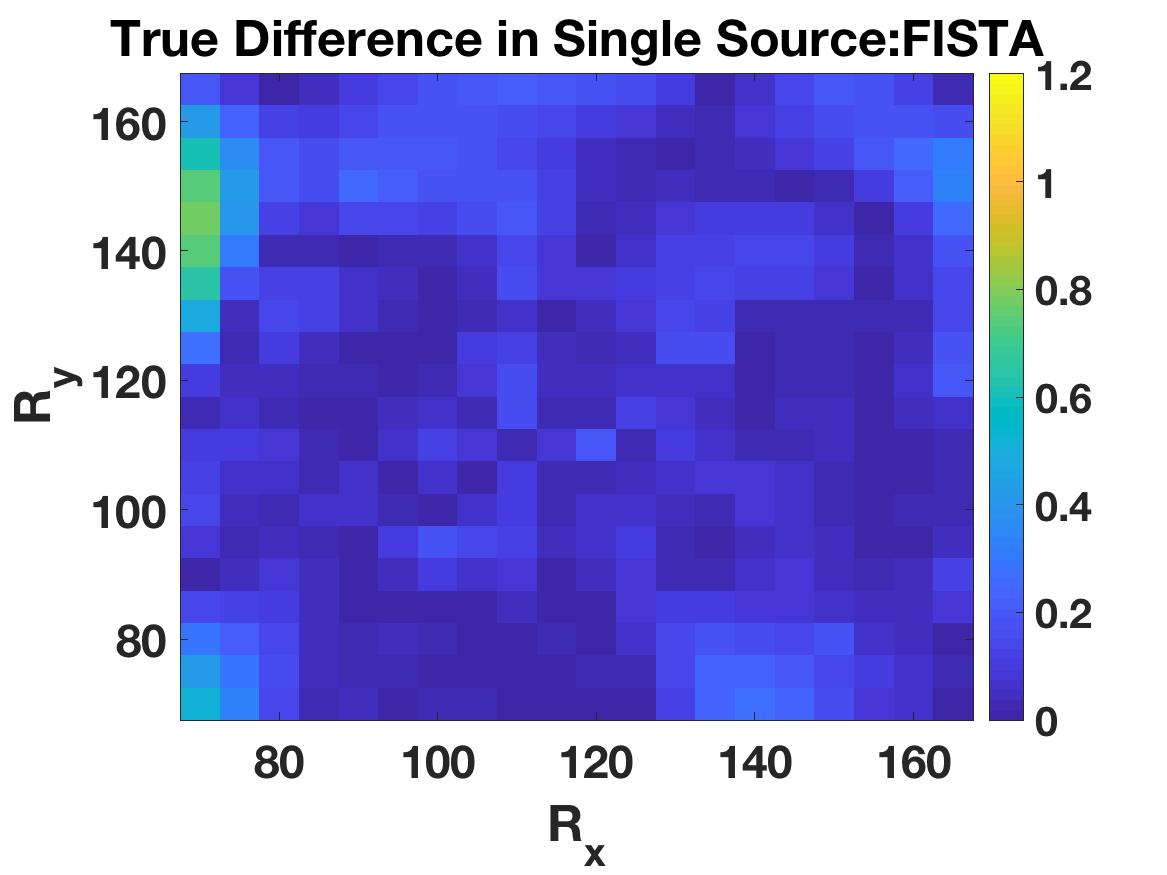}}
\subfigure[L-BFGS. ]{\label{fig:lbfgs_diff_single}\includegraphics[width=60mm]{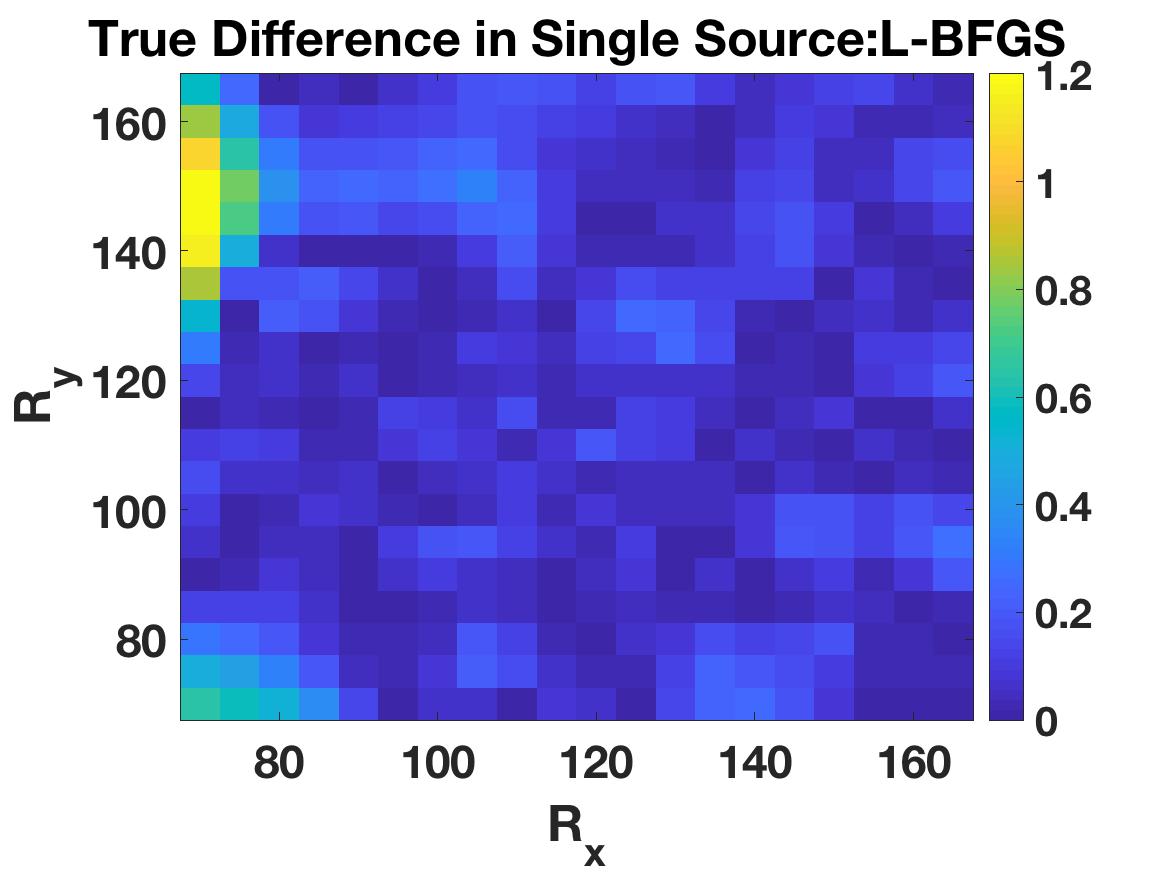}}
\subfigure[Low-Rank only.]{\label{fig:lr_diff_single}\includegraphics[width=60mm]{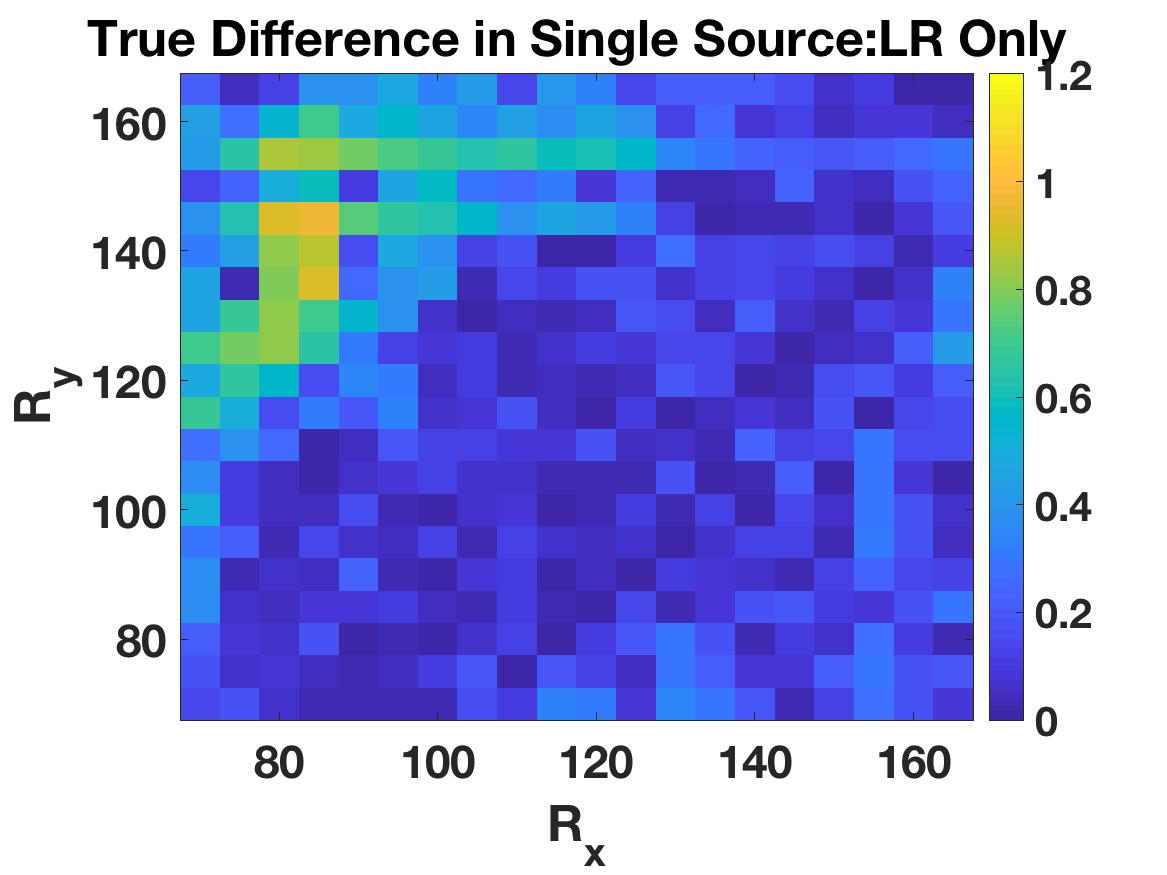}}
\subfigure[Smoothing only. ]{\label{fig:sm_diff_single}\includegraphics[width=60mm]{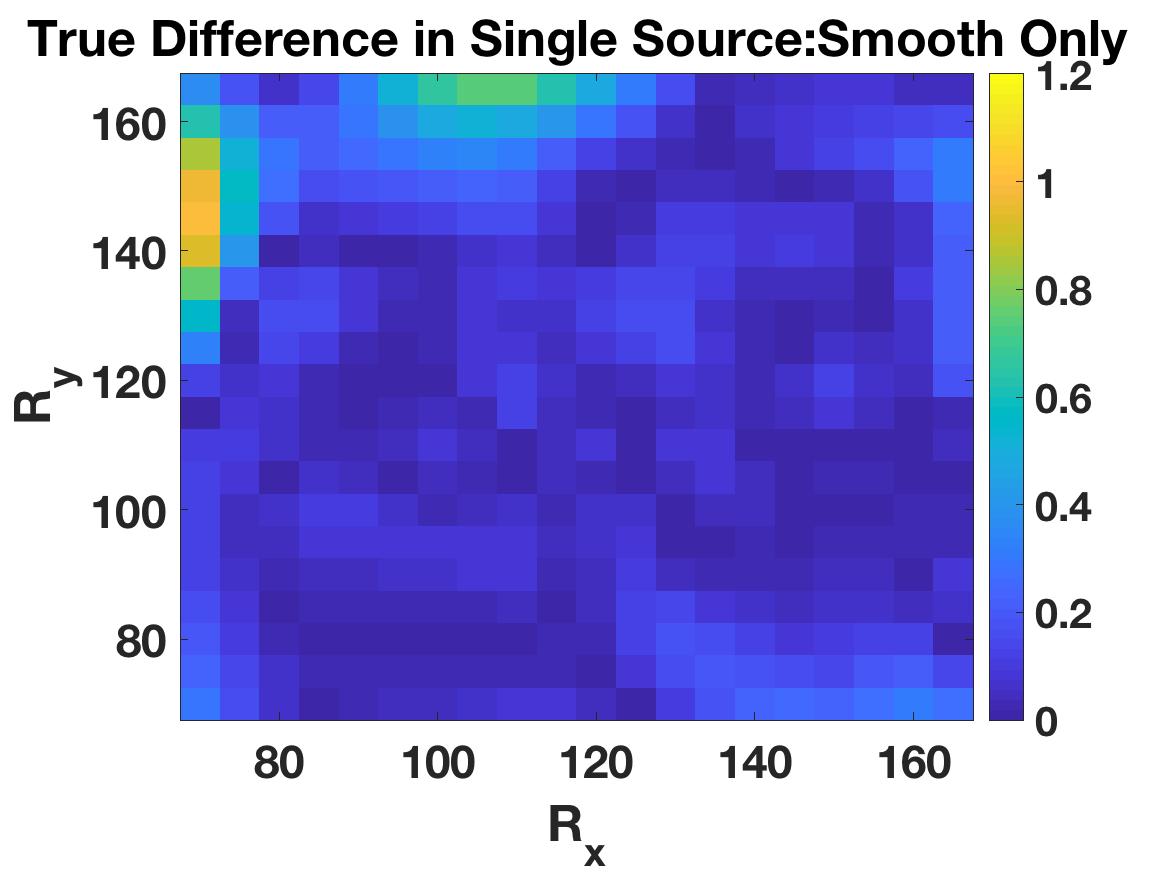}}
\caption{$|X - X_{true}|$ for all algorithms in a single source.}
\label{fig:algcomp_diff_single}
\end{figure}

With both smoothness and low-rank regularization, fitting observed data inexactly yields a better interpolation for missing data. 
Figure \ref{fig:algcomp_all} shows that variable relaxation inexact data fit (subfigure \ref{fig:vr_int}) 
has fewer contrasts overall when compared to the fits obtained using exact fit and competing algorithms. 
Focusing on a single source in Figure \ref{fig:algcomp_single}, we can see that the inequality constraint produces a smoother image for each particular source. 
While each algorithm effectively captures high energy area in the northwest corner of the plot, 
most algorithms overestimate the amount of energy that is actually present; variable relaxation with $\sigma>0$ gets the closest values. 
All algorithms tend to smooth out the observed entries, and tend to be less accurate the fewer entries there are in a designated space, notably around the corners of the receiver grid. Generally speaking, the combination of smoothness and low-rank works much better for interpolating the interior station grid rather than extrapolation near the edges.
%Notice also that the true data can be highly corrupted from the original data, as the magnitudes in each interpolation result and the ``true'' value in Figure \ref{fig:model_true_alg_single2} can be quite different. 
Figures \ref{fig:algcomp_diff_all} (full tensor) and \ref{fig:algcomp_diff_single} (single source) depict the absolute values of the results for each algorithm and the true value. Figure \ref{fig:algcomp_diff_all} shows that some sources are estimated very poorly, with the main error contributions coming further away from the mountain. 

Zooming in to a particular source, Figure \ref{fig:algcomp_diff_single} shows how much each algorithm overestimates the high-energy in the northwest corner of the map, primarily because the observed data in that corner is also very large. 
Certain artifacts of the interpolation schemes are seen more vividly in the difference plots. 
For instance, slightly northeast of the center, there is a higher energy region where very few data points are collected. 
This is very pronounced in the L-BFGS case (Figure \ref{fig:lbfgs_diff_single}) but less pronounced in the variable relaxation case (Figure \ref{fig:vr_diff_single}).  
Overall, the amount of error present in the $\sigma>0$ variable relaxation case is lower than all the other schemes. \\


\begin{figure}[H]
\centering     %%% not \center
\subfigure[Data misfit decay. ]{\label{fig:obj_total}\includegraphics[width=60mm]{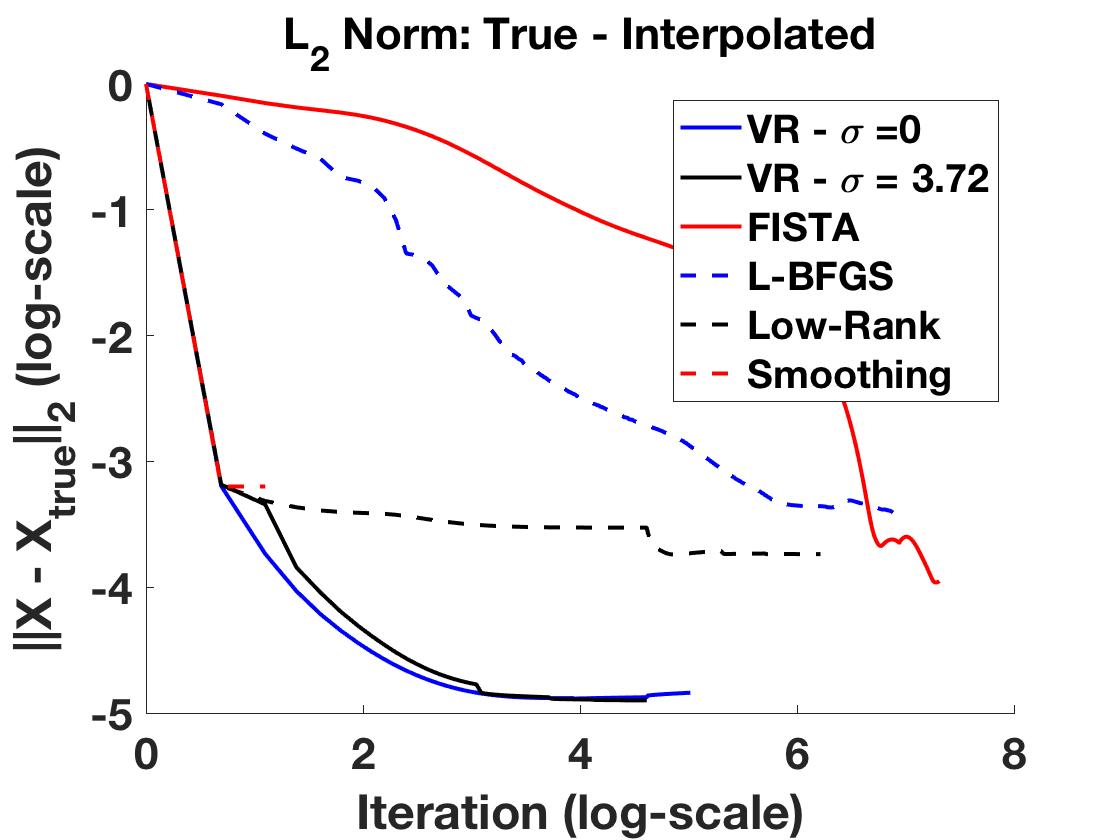}}
\subfigure[SVD decay for interpolated matrices.]{\label{fig:svd_total}\includegraphics[width=60mm]{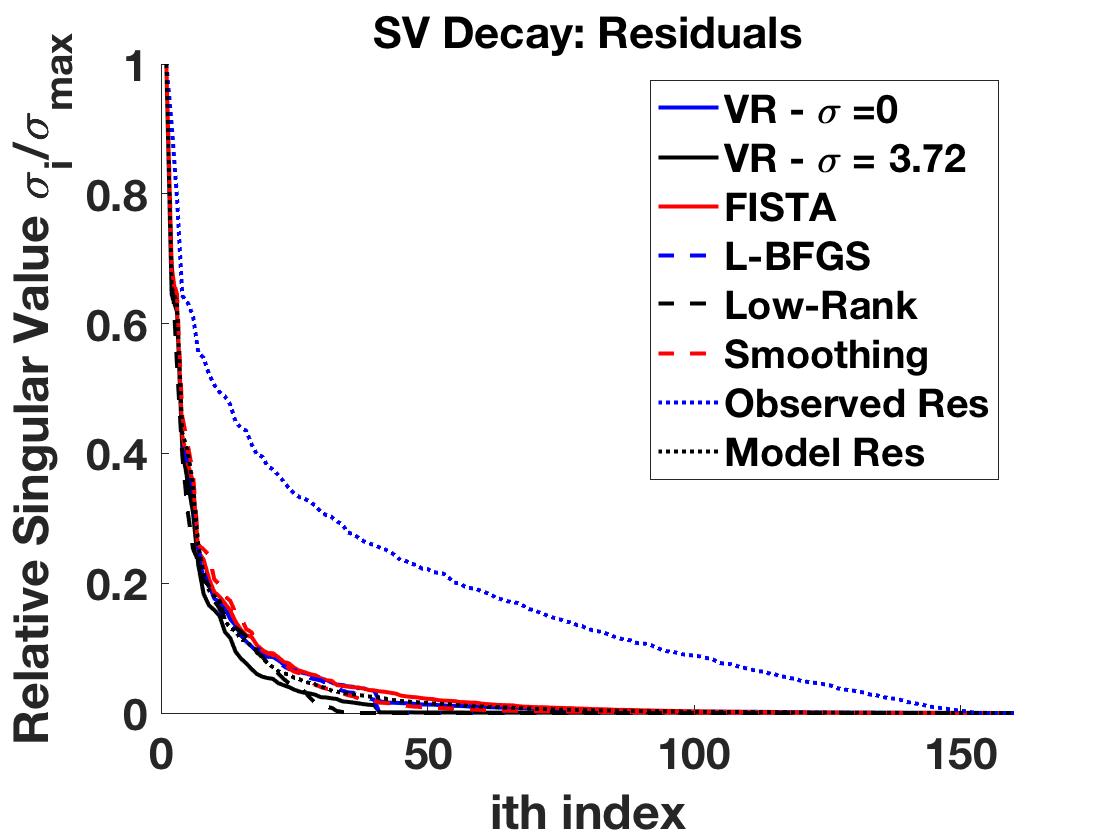}}
\caption{Convergence information for different algorithms.}
\label{fig:conv}
\end{figure}

Finally, we turn to singular value decay and convergence history. The singular value decay is shown in Figure \ref{fig:svd_total}, and a log-log convergence plot 
appears in Figure \ref{fig:obj_total}. With the exception of the smoothing-only implementation, all algorithms match the SVD decay well. In Figure \ref{fig:obj_total}, L-BFGS and FISTA have similar convergence histories while L-BFGS is relatively slow compared to the FISTA (and indeed all other algorithms) in compute time. 
%All of the variable relaxation/smoothing algorithms approach minima relatively quickly, albeit not necessarily the global minima since the problem is not convex. 
The proposed approach converges faster than competing methods, and also matches the SVD decay of the ground truth datasets.

% This is due to the nature of the gradient computation of $L$ and $R$. To find the gradient with respect to $L$ and $R$, we can rewrite the L-BFGS formulation in Equation ~\eqref{eq:lbfgs} using the vec-kron property. The gradient of $L$ can be computed by first rewriting $F$ using the classic kronecker identity $\mbox{vec}(AXB) = (B^T \otimes A) \mbox{vec}(X)$:
% \begin{align*}
% 	F(L, R) & = \tau\frac12\|L\|_F^2 + \tau\frac12\|R\|_F^2 + \frac{1}{2\gamma}\|\cL\vecp(LR^T)\|_2^2 + \|\cA\vecp(LR^T) - b\|_2^2 \\
% 			& = \tau\frac12\|\vecp(L)\|_F^2 + \tau\frac12\|\vecp(R)\|_F^2 + \frac{1}{2\gamma}\|\cL(R\otimes I)(\vecp(L)\|_2^2 + \|\cA(R\otimes I)\vecp(L) - b\|_2^2
% \end{align*}
%which depends on the last iteration of $R^T$. The gradient of $R$ is similar: 
% The derivatives are now straightfoward: 
% \begin{align*}
% 	\nabla_{\vecp(L)}F(L, R)& = \lt[\tau I + \frac{1}{\gamma}\lt(\cL(R\otimes I)^T\cL(R\otimes I)\rt)+ \lt((\cA(I\otimes L))^T(\cA(I\otimes L))\rt) \rt]\vecp(L) - (\cA(I\otimes L))^Tb.\\
% 	\nabla_{\vecp(R^T)}F(L, R)& = \lt[\tau I + \frac{1}{\gamma}\lt((I\otimes L^T)\cL^T\cL(I\otimes L)\vecp(R^T) \rt) + (I\otimes L^T)\cA^T\cA(I\otimes L)\vecp(R^T) \rt].
% \end{align*}
% \textcolor{red}{
% Of course, we can pre-compute $\cL^T\cL$ and $\cA^T\cA$, but we also must compute the rather costly kronecker product at each iteration, which slows down computational time. 
% }

\section{Conclusions and Future Directions}
Travel time tomography suffers from data collection constraints, reducing model resolution. 
We proposed an interpolation scheme that combines both local smoothness and low-rank information 
from a carefully chosen tessellation of the data 
to estimate residual values at prospective stations. 
To implement the scheme,  we developed a new relaxation approach that is flexible enough to allow multiple regularizers, 
and efficient in practice. We used it to estimate data from missing stations with a relatively high degree of accuracy (measured against a synthetic ground-truth dataset), 
in the presence of observation noise in available data. 
The algorithm is more flexible than available alternatives, and in particular can fit available data to a prescribed error level. 
The new approach is competitive with standard alternatives,
and offers new functionality to interpolate with both local and global structure over data fitting constraints. It is important to note that none of the methods do well estimating where there are very few data stations. This can be observed in Figure \ref{fig:algcomp_diff_single}, where the left most row between $R_y=140$ and $R_y = 160$. The proximity of at least one station increases the accuracy of our scheme.

Our next steps are to interpolate noisy data on a non-uniform grid and see if we can improve resolution in the seismic velocity inverse problem. 
In the experiments we presented, we used a synthetic dataset as `ground truth' to evaluate interpolation algorithms. 
The ultimate goal is to interpolate `missing' stations not present in the dataset, and see whether we can improve resolution 
of challenging regions within the upper crust near Mount St. Helens. 

%A future goal would be to develop a method of projecting data points onto the grid that does not drastically increase the error of interpolating new receivers, which can also be at off-grid locations. In addition, most travel times come with uncertainties observed at the time of data collection, and currently this method has no way of producing uncertainties for interpolated points. This may be done with a bootstrapping technique via sampling a subset of the already scarce data, or by some metric of distance to observed data. Finally, the ultimate goal of this technique would be to generate quality seismogram estimations efficiently for use in an inverse problem, with the hope of achieving increased accuracy for the end result. 

\section{Acknowledgements}
R. Baraldi acknowledges support from the Department of Energy Computational Science Graduate Fellowship, which is provided under grant number DE-FG02-97ER25308. 
C. Ulberg and K. Creager were supported by National Science Foundation grant number EAR-1358512.
The work of A. Aravkin was supported by the Washington Research Foundation Data Science Professorship.

\section{References}
\bibliographystyle{abbrv}  
\bibliography{vr_ref}

\begin{thebibliography}{10}

\bibitem{aravkin2014fast}
A.~Aravkin, R.~Kumar, H.~Mansour, B.~Recht, and F.~J. Herrmann.
\newblock Fast methods for denoising matrix completion formulations, with
  applications to robust seismic data interpolation.
\newblock {\em SIAM Journal on Scientific Computing}, 36(5):237--266, 2014.

\bibitem{aravkin2016level}
A.~Y. Aravkin, J.~V. Burke, D.~Drusvyatskiy, M.~P. Friedlander, and S.~Roy.
\newblock Level-set methods for convex optimization.
\newblock {\em arXiv preprint arXiv:1602.01506}, 2016.

\bibitem{beck2009FISTA}
A.~Beck and M.~Teboulle.
\newblock A fast iterative shrinkage-thresholding algorithm for linear inverse
  problems.
\newblock {\em SIAM Journal on Imaging Sciences}, 2(1):183--202, 2009.

\bibitem{candes2010Matcompnoise}
E.~J. Cand{\`e}s and Y.~Plan.
\newblock Matrix completion with noise.
\newblock {\em Proceedings of the IEEE}, 98(6):925--936, June 2010.

\bibitem{Candes2009matcompconvex}
E.~J. Cand{\`e}s and B.~Recht.
\newblock Exact matrix completion via convex optimization.
\newblock {\em Foundations of Computational Mathematics}, 9(6):717--723, Apr
  2009.

\bibitem{candes2006near}
E.~J. Cand{\`e}s and T.~Tao.
\newblock Near-optimal signal recovery from random projections: universal
  encoding strategies.
\newblock {\em IEEE Transactions on Information Theory}, 52(12):5406--5425,
  2006.

\bibitem{splittingschemes}
D.~Davis and W.~Yin.
\newblock Convergence rate analysis of several splitting schemes.
\newblock In {\em Splitting Methods in Communication, Imaging, Science, and
  Engineering}, pages 115--163. Springer, 2016.

\bibitem{hennenfent2008simply}
G.~Hennenfent and F.~J. Herrmann.
\newblock Simply denoise: wavefield reconstruction via jittered undersampling.
\newblock {\em Geophysics}, 73(3):19--28, 2008.

\bibitem{herrmann2008non}
F.~J. Herrmann and G.~Hennenfent.
\newblock Non-parametric seismic data recovery with curvelet frames.
\newblock {\em Geophysical Journal International}, 173(1):233--248, 2008.

\bibitem{hole19953dfdreflection}
J.~A. Hole and B.~C. Zelt.
\newblock 3-{D} finite-difference reflection traveltimes.
\newblock {\em Geophysical Journal International}, 121(2):427--434, 1995.

\bibitem{kiser2016vpvs}
E.~Kiser, A.~Levander, C.~Zelt, B.~Schmandt, and S.~Hansen.
\newblock Focusing of melt near the top of the {M}ount {S}t. {H}elens ({USA})
  magma reservoir and its relationship to major volcanic eruptions.
\newblock {\em Geology}, 2018.

\bibitem{2017AGUkiser}
E.~{Kiser}, A.~{Levander}, C.~A. {Zelt}, I.~{Palomeras}, K.~{Creager}, C.~W.
  {Ulberg}, B.~{Schmandt}, S.~M. {Hansen}, S.~H. {Harder}, G.~A. {Abers}, and
  K.~{Crosbie}.
\newblock Three-dimensional velocity models of the {M}ount {S}t. {H}elens
  magmatic system using the i{MUSH} active-source data set.
\newblock In {\em AGU Fall Meeting Abstracts}, Dec 2017.

\bibitem{kumar2015efficient}
R.~Kumar, C.~Da~Silva, O.~Akalin, A.~Y. Aravkin, H.~Mansour, B.~Recht, and
  F.~J. Herrmann.
\newblock Efficient matrix completion for seismic data reconstruction.
\newblock {\em Geophysics}, 80(5):97--114, 2015.

\bibitem{landro2002uncertainties}
M.~Landr{\o}.
\newblock Uncertainties in quantitative time-lapse seismic analysis.
\newblock {\em Geophysical Prospecting}, 50(5):527--538, 2002.

\bibitem{landro2008effect}
M.~Landr{\o}.
\newblock The effect of noise generated by previous shots on seismic reflection
  data.
\newblock {\em Geophysics}, 73(3):9--17, 2008.

\bibitem{crosson1989helens}
J.~Lees and R.~Crosson.
\newblock Tomographic inversion for three‐dimensional velocity structure at
  {M}ount {S}t. {H}elens using earthquake data.
\newblock {\em Journal of Geophysical Research: Solid Earth},
  94(B5):5716--5728, 1989.

\bibitem{pavlis1983vstruct}
S.~D. Malone and G.~L. Pavlis.
\newblock Velocity structure and relocation of earthquakes at {M}ount {S}t.
  {H}elens.
\newblock In {\em EosTrans AGU}, volume~64, page 895, 1983.

\bibitem{martins2001connection}
J.~Martins, P.~Sturdza, and J.~Alonso.
\newblock The connection between the complex-step derivative approximation and
  algorithmic differentiation.
\newblock In {\em 39th Aerospace Sciences Meeting and Exhibit}, page 921, 2001.

\bibitem{moran1999rainier}
S.~C. Moran, J.~M. Lees, and S.~D. Malone.
\newblock $p$-wave crustal velocity structure in the greater {M}ount {R}ainier
  area from local earthquake tomography.
\newblock {\em Journal of Geophysical Research: Solid Earth},
  104(B5):10775--10786, 1999.

\bibitem{NocedalWright.optim.2000}
J.~{Nocedal} and S.~J. {Wright}.
\newblock {\em Numerical Optimization}.
\newblock Springer, 2000.

\bibitem{oneto2016tikhonov}
L.~Oneto, S.~Ridella, and D.~Anguita.
\newblock Tikhonov, {I}vanov and {M}orozov regularization for support vector
  machine learning.
\newblock {\em Machine Learning}, 103(1):103--136, 2016.

\bibitem{phillips1991comparison}
W.~S. Phillips and M.~C. Fehler.
\newblock Traveltime tomography: A comparison of popular methods.
\newblock {\em Geophysics}, 56(10):1639--1649, 1991.

\bibitem{Recht2010guaranteed}
B.~Recht, M.~Fazel, and P.~A. Parrilo.
\newblock Guaranteed minimum-rank solutions of linear matrix equations via
  nuclear norm minimization.
\newblock {\em SIAM Rev.}, 52(3):471--501, Aug 2010.

\bibitem{sacchi}
M.~D. Sacchi, T.~J. Ulrych, and C.~J. Walker.
\newblock Interpolation and extrapolation using a high-resolution {D}iscrete
  {F}ourier transform.
\newblock {\em IEEE Transactions on Signal Processing}, 46(1):31--38, Jan 1998.

\bibitem{tseng2001convergence}
P.~Tseng.
\newblock Convergence of a block coordinate descent method for
  nondifferentiable minimization.
\newblock {\em Journal of optimization theory and applications},
  109(3):475--494, 2001.

\bibitem{2017AGUFMcarl}
C.~W. {Ulberg}, K.~{Creager}, S.~C. {Moran}, G.~A. {Abers}, K.~{Crosbie}, R.~S.
  {Crosson}, R.~P. {Denlinger}, W.~A. {Thelen}, E.~{Kiser}, A.~{Levander}, and
  O.~{Bachmann}.
\newblock Imaging seismic zones and magma beneath {M}ount {S}t. {H}elens with
  the i{MUSH} {B}roadband {A}rray.
\newblock In {\em AGU Fall Meeting Abstracts}. Submitted, Dec 2017.

\bibitem{vidale19903dfd}
J.~E. Vidale.
\newblock Finite-difference calculation of traveltimes in three dimensions.
\newblock {\em Geophysics}, 55(5):521--526, 1990.

\bibitem{WAITE2009113}
G.~P. Waite and S.~C. Moran.
\newblock ${V}_p$ structure of {M}ount {S}t. {H}elens, {W}ashington, {USA},
  imaged with local earthquake tomography.
\newblock {\em Journal of Volcanology and Geothermal Research},
  182(1):113--122, 2009.

\bibitem{zheng2018fast}
P.~{Zheng} and A.~{Aravkin}.
\newblock Fast methods for nonsmooth nonconvex minimization.
\newblock {\em arXiv preprint arXiv:1802.02654}, 2018.

\end{thebibliography}

\end{document}